\documentclass[final,leqno]{article}

\usepackage{amsfonts,amssymb,amsmath,mathtools,color,graphics,mathdots}
\usepackage{amstext}
\usepackage{amscd}
\usepackage{amsthm}
\usepackage{graphicx}
\usepackage{longtable}
\usepackage{calc}
\usepackage{bbding}
\usepackage[algoruled]{algorithm2e}

\usepackage{booktabs}

\usepackage{tikz}

\newtheorem{theorem}{Theorem}
\newtheorem{lemma}[theorem]{Lemma}

\DeclarePairedDelimiter{\abs}{\lvert}{\rvert}
\DeclarePairedDelimiter{\norm}{\lVert}{\rVert}

\newsavebox{\savepar}

\newcommand{\wt}{\widetilde}

\newcommand{\iunit}{{\mathfrak i}}

\newcommand{\la}{\lambda}

\newcommand{\twoone}[2]{\left[\begin{array}{c} #1\\ #2 \end{array}\right]}
\newcommand{\twotwo}[4]{\left[\begin{array}{cc} #1 & #2 \\ #3 & #4 \end{array}\right]}

\DeclareMathOperator{\opvec}{vec}

\usepackage{todonotes}

\newcommand{\C}{\mathbb{C}}

\title{Invariant subspaces of $T$-palindromic pencils and algebraic $T$-Riccati equations}

\author{Bruno Iannazzo\thanks{Dipartimento di Matematica e Informatica, University of Perugia, Italy, 
{\tt bruno.iannazzo@unipg.it}, 
corresponding author
}, 
Beatrice Meini\thanks{Dipartimento di Matematica, University of Pisa, Italy, 
{\tt beatrice.meini@unipi.it}
}, 
Federico Poloni\thanks{Dipartimento di Informatica, University of Pisa, Italy, 
{\tt federico.poloni@unipi.it}
}
}

\usepackage[notcite,notref]{showkeys}

\begin{document}

\maketitle

\begin{abstract}
By exploiting the connection between solving algebraic $\top$-Riccati equations and computing certain deflating subspaces of $\top$-palindromic matrix pencils, we obtain theoretical and computational results on both problems. Theoretically, we introduce conditions to avoid the presence of modulus-one eigenvalues in a $\top$-palindromic matrix pencil and conditions for the existence of solutions of a $\top$-Riccati equation. Computationally, we improve the palindromic QZ algorithm with a new ordering procedure and introduce new algorithms for computing a deflating subspace of the $\top$-palindromic pencil, based on quadraticizations of the pencil or on an integral representation of the orthogonal projector on the sought deflating subspace.\\

\end{abstract}

\section{Introduction}

We consider the Nonsymmetric Algebraic $\top$-Riccati equation ($\top$-NARE)
\begin{equation}\label{eq:tnare}
  \mathcal R(X)=0,\qquad \mathcal R(X):=DX + X^\top A-X^\top BX+C,
\end{equation}
where $A,B,C,D\in\mathbb{R}^{n\times n}$ are given and the $n\times n$ matrix $X$ is the unknown, while the superscript $\top$ denotes transposition. 

The nonlinear equation \eqref{eq:tnare} has been studied in \cite{ikr} where the author characterizes the solution in terms of the neutral space of a suitable matrix; in \cite{vor} where some sufficient conditions are given for the existence of the solution in terms of nonsingularity of the matrix coefficients; in \cite{yuan,yuan1} where the matrix coefficients are complex and satisfy suitable Hermitianity properties. Concerning numerical methods, in \cite{vor2} a numerical method for the computation of the solution is proposed under suitable symmetry assumptions, in \cite{bp20} a fixed point iteration and Newton's method are proposed, while in \cite{miy} interval arithmetic is used to compute enclosures of a solution. 
The case with $A=D=0$ has been treated in \cite{borobia,borobia2}. Examples of applications with $A=D=0$ are encountered in image restoration problems \cite{edb}, in inverse eigenvalue problems \cite{aishima}, and in the numerical solution of a cubic matrix equation arising in conservative dynamics~\cite{bv}.

As in \cite{bimp}, we associate with  equation \eqref{eq:tnare} the $2n\times 2n$ matrix
\begin{equation}\label{eq:M}
\mathcal M=\begin{bmatrix} C & D\\A & -B\end{bmatrix},
\end{equation}
and the {\em $\top$-palindromic }pencil
\begin{equation}\label{eq:phi}
\varphi(z)=\mathcal M + z \mathcal M^\top.
\end{equation}
We may easily observe that, if $X$ is such that $\mathcal R(X)=0$, then 
\begin{equation}\label{eq:phial}
\varphi(z)\begin{bmatrix} I \\ X \end{bmatrix}	=\begin{bmatrix} -X^\top\\I\end{bmatrix}
	\alpha(z),~~\alpha(z)=A-BX+z(D^\top-B^\top X).
\end{equation}
In \cite{bimp} it is shown that, under suitable assumptions on the pencil $\varphi(z)$, a converse result is true. 
This characterization leads to new numerical methods for approximating the solution of \eqref{eq:tnare}, based on the computation of an invariant subspace of the pencil $\varphi(z)$. In particular, algorithms based on the (palindromic) QZ method and on the Doubling Algorithm are proposed in \cite{bimp}.

In this paper we rely on the relationship between solutions of \eqref{eq:tnare} and deflating subspaces of the pencil $\varphi(z)$ to obtain some advances concerning theoretical properties and develop new numerical methods for both solving $\top$-NAREs 
and finding deflating subspaces of palindromic pencils.

Most of the results and numerical methods
require the assumption that $\varphi(z)$ is nonsingular if $|z|=1$. However, this property is not always easy to verify. To this regard, we introduce some sufficient conditions under which $\det\varphi(z)\ne 0$ for $|z|=1$, which resemble the well-known conditions ensuring the absence of critical eigenvalues in algebraic Riccati equations.
Concerning the existence of solutions, we provide new sufficient conditions, expressed in terms of the quadratic matrix coefficient $B$. A necessary condition is also introduced, which resembles a similar property for continuous-time algebraic Riccati equations.

On the computational side, we introduce an improvement of the numerical method proposed in \cite{bimp} based on palindromic QZ  algorithm, and we introduce some new algorithms.
The improvement in the palindromic QZ algorithm consists in applying a more efficient swapping strategy for the eigenvalues in the antritriangular form of a $\top$-palindromic pencil. Indeed, such swapping is needed to select the deflating subspace corresponding to the eigenvalues lying in the open unit disk. The proposed modification reduces the computational cost of the swapping procedure with respect to the strategy followed in \cite{bimp}. In particular, the number of arithmetic operations is reduced from $\mathcal{O}(n^4)$ to $\mathcal{O}(n^3)$ in the worst case. An additional benefit is that the numerical stability is improved.
 
Two new algorithms consist in transforming the pencil $\mathcal \varphi(z)$ into a quadratic matrix polynomial  and in the application of Cyclic Reduction \cite{bm:cr}. The two algorithms differ in the way the quadratic matrix polynomials are constructed. In the first approach (\textbf{CR1}), the quadratic polynomial is obtained by adding  $n$ eigenvalues at 0 and $n$ eigenvalues at infinity to the pencil $\varphi(z)$. In the second approach (\textbf{CR2}) the matrix polynomial is $\top$-palindromic. Both algorithms have quadratic convergence and a computational cost of $\mathcal{O}(n^3)$ arithmetic operations per step. The overhead constant in the second algorithm is about half as large as the constant in the first algorithm, but the second algorithm might suffer from numerical instability in ill-conditioned problems.
 
A third new algorithm (denoted by \textbf{Int}) relies on a contour integral representation of the orthogonal projector on the deflating subspace of $\varphi(z)$ corresponding to the eigenvalues in the open unit disk, in a strategy similar to that of FEAST eigensolvers \cite{feast}. The integral is approximated by the trapezoidal scheme, by evaluating the integrand at the complex roots of 1. The resulting algorithm reveals connections with DFT and with doubling algorithms for computing the stable deflating subspace of a pencil. The overall computational cost is $\mathcal{O}(Kn^3)$, where $K$ in the number of nodes used in the integral approximation.
 
The proposed algorithms have been tested on a variety of problems, with different peculiarities. The experiments show that the palindromic QZ method is the most accurate algorithm. Indeed, even if applied to severely ill-conditioned problems, the accuracy of the approximation remains of the order of the machine precision. The algorithm based on the integral representation \textbf{Int} and algorithm \textbf{CR2} are methods which exploit the symmetry of the problem. However, algorithm \textbf{Int} is more accurate than \textbf{CR2}. On the other hand, the latter algorithm is faster than the former. The algorithm \textbf{CR1} performs similarly to the Doubling Algorithm proposed in \cite{bimp}.

The paper is organized as follows. In Section~\ref{sec:prel} we recall some definition and properties related to matrix pencils. Sections~\ref{sec:cri}
and \ref{sec:ex} concern theoretical properties of the $\top$-NARE \eqref{eq:tnare}: in the former sufficient conditions are presented for the nonsingularity of $\varphi(z)$ on the unit circle, in the latter sufficient and necessary conditions are introduced for the existence of a solution. In Section~\ref{sec:cl} we recall classical algorithms, for comparisons with the newly introduced ones. In Section~\ref{sec:sw} the new strategy for swapping the eigenvalues is described. Section~\ref{sec:qua} deals with quadraticizations of the pencil $\varphi(z)$ and with the application of Cyclic Reduction for the solution of the $\top$-NARE \eqref{eq:tnare}. In Section~\ref{sec:ir} the method based on the integral representation is proposed. Numerical experiments are presented in Section~\ref{sec:exp}, while conclusions are drawn in Section~\ref{sec:con}.

\section{Preliminaries}\label{sec:prel}

We recall some definitions and properties regarding matrix pencils, which will be used in the following. For a thorough treatise on this topic, we refer to~\cite{glr}.

A matrix pencil $p(z)=A+zB$, with $A,B\in\mathbb C^{n\times n}$, such that $\det(p(z))$ is not the zero polynomial is said to be regular. When the pencil $p(z)$ is regular, the zeros of $\det(p(z))$ are its finite eigenvalues, while infinity has multiplicity $d\ge 1$ as eigenvalue of $p(z)$ when $\det(p(z))$ has degree $n-d$; the spectrum of $p(z)$ is the set of its eigenvalues. Observe that for a regular pencil, $0$ [resp. $\infty$] is an eigenvalue if and only if $A$ [resp. $B$] is singular.

A subspace $\mathcal{V}\subset \mathbb{C}^n$ of dimension $\ell$ is a deflating subspace of the regular matrix pencil $p(z)=A+zB$ if there exists a subspace $\mathcal{W}\subset \mathbb{C}^n$ of dimension $\ell$ such that $A\mathcal{V} \subset \mathcal{W}$ and $B\mathcal{V} \subset \mathcal{W}$. We call the subspace $\mathcal W$ image of the deflating subspace $\mathcal V$. If the columns of $V\in\mathbb{C}^{n\times \ell}$ and $W\in\mathbb{C}^{n\times \ell}$ provide a basis for $\mathcal{V}$ and $\mathcal{W}$, respectively, then there exist $A',B'\in\mathbb{C}^{\ell\times \ell}$ such that $AV=WA'$ and $BV=WB'$. The $\ell\times \ell$ pencil $A'+zB'$ is regular and its eigenvalues (which are a subset of the eigenvalues of $p(z)$) are said to be the spectrum or the eigenvalues of $p(z)$ associated with the deflating subspace $\mathcal{V}$. We say also that $\mathcal V$ is the deflating subspace associated with the eigenvalues of $A'+zB'$.

A graph deflating subspace of $p(z)$ is a deflating subspace $\mathcal{V}$ such that any matrix $V\in\mathbb C^{n\times \ell}$ whose columns are a basis for $\mathcal{V}$ has its leading $\ell\times \ell$ submatrix invertible. Equivalently, there exists $X\in\C^{(n-\ell)\times \ell}$ such that the columns of the matrix $\left[\begin{smallmatrix} I_\ell\\X\end{smallmatrix}\right]$ are a basis of $\mathcal{V}$.

For a matrix $M$, the notation $M^\top$ denotes its transpose, and for a matrix $M$ with complex entries the notation $M^*$ denotes its transpose conjugate. We write $M \succ 0$ (resp. $M \succeq 0$) to denote that a Hermitian matrix $M$ is positive definite (resp. positive semidefinite).

We denote by $\mathbb{T}$ the unit circle, i.e., $\mathbb{T}=\{z\in\mathbb C\; : \: |z|=1\}$. 

A pencil $p(z)=A+zB$ is said to be $\top$-palindromic if $B=A^\top$. If $\lambda$ is an eigenvalue of a regular $\top$-palindromic pencil, then so is $1/\lambda$, and the two have the same algebraic and geometric multiplicities; this property holds also for $0$ and $\infty$, with the conventions $1/0=\infty$ and $1/\infty=0$.
A $\top$-palindromic pencil $p(z)$ is said to be non-critical if $\det(p(z))\ne 0$ for any $z\in\mathbb{T}$, critical otherwise. 
A non-critical $\top$-palindromic pencil $p(z)$ of size $2n$  has $n$ eigenvalues inside the unit disk, associated with a unique deflating subspace, said to be the stable deflating subspace.

Two pencils $p(z)=A+zB$ and $\wt p(z)=\wt A+z\wt B$ are said to be right equivalent if there exists an invertible matrix $M$ such that $\wt A=MA$ and $\wt B=MB$. Two right equivalent pencils share the same eigenvalues and deflating subspaces. 

A set $\mathcal{S}\subset \mathbb{C}\cup\{\infty\}$ is said to be reciprocal-free when it does not contain a pair of elements in the form $(\lambda, 1/\lambda)$ for some $\lambda \in \mathbb{C} \cup \{\infty\}$.

A matrix $M\in\C^{n\times n}$ is antitriangular if $(M)_{ij}=0$ when $i+j\le n$. The antidiagonal of $M$ is the set of elements with $i+j=n+1$.
A pencil $A+zB$ is said to be in antitriangular form if both $A$ and $B$ are antitriangular matrices. An antitriangular pencil is regular if and only if $A_{i,n+1-i}$ and $B_{i,n+1-i}$ are not both zero for all $i=1,\ldots,n$. The eigenvalues of a regular antitriangular pencil can be read from its antidiagonal as the solutions of
the scalar equations $A_{i,n+1-i}+zB_{i,n+1-i}=0$ (including infinity if $B_{i,n+1-i}=0$).

Finally, we recall the following result \cite[Theorem 2]{bimp}.

\begin{theorem}\label{thm:thm2}
Assume that the pencil $\varphi(z)$ of \eqref{eq:phi} is regular.
If $X$ is an $n\times n$ matrix such that the columns of the matrix 
$\left[\begin{smallmatrix} I_n \\X\end{smallmatrix}\right]$
span a
deflating subspace of $\varphi(z)$ associated with a reciprocal-free set of eigenvalues, then
$X$ is a solution to~\eqref{eq:tnare}.
\end{theorem}

From \eqref{eq:phial}, if $X$ is a solution of \eqref{eq:tnare}
then the columns the matrix 
$\left[\begin{smallmatrix} I_n \\X\end{smallmatrix}\right]$
span a
deflating subspace of $\varphi(z)$.
A solution $X$ to \eqref{eq:tnare} is said to be stabilizing if such
deflating subspace is associated with eigenvalues of $\varphi(z)$ lying in the open unit disk.

\section{Critical eigenvalues}\label{sec:cri}

In this section, we discuss conditions under which $\varphi(z)$, defined in \eqref{eq:phi}, is non-critical, and hence it has $n$ eigenvalues inside the unit circle an $n$ outside.

We start from a classical lemma on saddle-point matrices. 
\begin{lemma} \label{lem:saddle}
Let
\[
\mathcal{A} = \begin{bmatrix}
\mathcal{C} & \mathcal{E}\\
\mathcal{E}^* & -\mathcal{B}
\end{bmatrix} \in \mathbb{C}^{(n+m)\times (n+m)}.
\]
If one of the following conditions holds
\begin{enumerate}
    \item $\mathcal{C} = \mathcal{C}^* \succ 0$ and $\mathcal{B} = \mathcal{B}^* \succ 0$;
    \item $\mathcal{C} = \mathcal{C}^* \succeq 0$, $\mathcal{B} = \mathcal{B}^* \succeq 0$ and $\begin{bmatrix}\mathcal{C} & \mathcal{E}\end{bmatrix}$ has full rank;
    \item $\mathcal{C} = \mathcal{C}^* \succeq 0$, $\mathcal{B} = \mathcal{B}^* \succeq 0$ and $\begin{bmatrix}\mathcal{E}^* & -\mathcal{B}\end{bmatrix}$ has full rank;
\end{enumerate}
then $\mathcal{A}$ is nonsingular.
\end{lemma}
\begin{proof}
Suppose otherwise that
\[
v = \begin{bmatrix}v_1 \\ v_2 \end{bmatrix} \neq 0 
\]
is in the kernel of $\mathcal{A}$.
Then, we can compute
\[
0 = \operatorname{Re} \left(\begin{bmatrix}v_1^* & -v_2^*\end{bmatrix} 
\begin{bmatrix}
\mathcal{C} & \mathcal{E}\\
\mathcal{E}^* & -\mathcal{B}
\end{bmatrix}
\begin{bmatrix}v_1 \\ v_2 \end{bmatrix}\right)
= \underbrace{v_1^*\mathcal{C}v_1}_{\geq 0} + \underbrace{v_2^* \mathcal{B} v_2}_{\geq 0} + \underbrace{\operatorname{Re}\left( v_1^*\mathcal{E}v_2 - v_2^*\mathcal{E}^* v_1\right)}_{=0}.
\]
By positive semidefiniteness we must have 
\begin{equation} \label{conditions}
\mathcal{C}v_1 = \mathcal{B} v_2 = 0.    
\end{equation}
If $\mathcal{B}$ and $\mathcal{C}$ are both positive definite, then we obtain $v_1=v_2=0$, which contradicts the assumption on $v$. Otherwise, plugging~\eqref{conditions} into $\mathcal{A}v = 0$ we obtain $\mathcal{E}v_2 = 0$ and $\mathcal{E}^*v_1 = 0$, i.e., $0 = v_1^*\begin{bmatrix}\mathcal{C} & \mathcal{E}\end{bmatrix} = v_2^*\begin{bmatrix}\mathcal{E}^* & -\mathcal{B}\end{bmatrix}$.
\end{proof}
Using this lemma, we can give sufficient conditions to ensure that the eigenvalues of $\varphi(z)$ do not lie on the unit circle. 

\begin{theorem}
Let $\varphi(z)$ be the $\top$-palindromic pencil defined in \eqref{eq:phi}, and let one of the following conditions holds:
	    \begin{enumerate}
	        \item $C=C^\top \succ 0$ and $B=B^\top \succ 0$; and moreover $D-A^\top$ has full rank;
	        \item $C=C^\top \succeq 0$ and $B=B^\top \succeq 0$; and moreover for each $y\in\mathbb{T}$ the matrix $\begin{bmatrix} y^{-1}C + yC^\top & y^{-1}D + yA^\top \end{bmatrix}$ has full rank;
	        \item $C=C^\top \succeq 0$ and $B=B^\top \succeq 0$; and moreover for each $y\in\mathbb{T}$ the matrix $\begin{bmatrix} y^{-1}A + yD^\top & -y^{-1}B - yB^\top \end{bmatrix}$ has full rank.
	    \end{enumerate}
	    Then, $\varphi(z)$ is regular and has no eigenvalues on the unit circle.
	\end{theorem}
	\begin{proof}
	    Suppose $\varphi(z)$ is singular for some $z_* \in \mathbb{T}$. Suppose for now $z_* \neq -1$. Then there exists $y$ with $\operatorname{Re}(y) > 0$ such that $y^2 = z_*$, and the matrix $\mathcal{A} = y^{-1} \mathcal{M} + y \mathcal{M}^\top = y^{-1}\varphi(z_*)$ is singular, too. Note that $\mathcal{A}$ is Hermitian, since $y^{-1} = \overline{y}$ when $y\in\mathbb{T}$, and that $y^{-1}C+yC^\top = 2 \operatorname{Re}(y) C$ is positive definite in case 1 and semidefinite in cases 2 and 3, and the same holds for $y^{-1}B+yB^\top$. The result then follows by applying Lemma~\ref{lem:saddle} to~$\mathcal{A}$. 
	    
	    In the case $z_*=-1$, we can still argue as above in cases 2 and 3, since $y^{-1}C+yC^\top = 0 \succeq 0$; in case 1, $\varphi(z_*)= \begin{bsmallmatrix}0 & D-A^\top \\ A-D^\top & 0\end{bsmallmatrix}$ is invertible because $D-A^\top$ has full rank.
	\end{proof}
	We note that Conditions~2 and~3 are formally similar to the classical conditions that ensure the absence of critical eigenvalues for the algebraic Riccati equations
	\[
	    XGX-A^\top X-XA -Q=0,
	\]
with $Q\succeq 0, G\succeq 0$ and $\begin{bmatrix}A & G\end{bmatrix}$  stabilizable~\cite[Theorem~7.2.8]{lancasterrodman-are}, and are indeed obtained by similar arguments.

\section{Existence of solutions}\label{sec:ex}

Let $\varphi(z)$ be regular and let $\mathcal{U}$ be a deflating subspace of $\varphi(z)$ associated with a reciprocal-free subset $\mathcal{S}$ of eigenvalues of $\alpha(z)$. 

\begin{theorem}[\protect{\cite[Theorem~3.3]{m4}}]
Let $\varphi(z)$ be a regular $\top$-palindromic pencil and let $\mathcal{U}$ be a deflating subspace of $\varphi(z)$
associated with a reciprocal-free subset $\mathcal{S}$ of eigenvalues of $\alpha(z)$ and with image $\mathcal{V}$. Then, $\mathcal{U}$ is isotropic, i.e., $u^\top \varphi(z) u = 0$ for each $u\in\mathcal{U}$.
\end{theorem}

We can write
\begin{equation}\label{eq:20}
\varphi(z)\begin{bmatrix} U_1\\U_2\end{bmatrix}
=\begin{bmatrix} V_1\\V_2\end{bmatrix} \wt \alpha(z),
\end{equation}
where the columns of $\left[\begin{smallmatrix} U_1\\U_2\end{smallmatrix}\right]$ and $\left[\begin{smallmatrix} V_1\\ V_2\end{smallmatrix}\right]$ span the subspace $\mathcal{U}$ and $\mathcal{V}$, respectively. Then, isotropy is equivalent to
\begin{equation}
    \begin{bmatrix} U_1^\top &U_2^\top\end{bmatrix} \varphi(z)\begin{bmatrix} U_1\\U_2\end{bmatrix} = 0.
\end{equation}

We look for conditions on the pencil $\varphi(z)$ of size $2n$ and the deflating subspace $\mathcal U$ of dimension $n$, under which $U_1$ is invertible and thus there exists a solution of the $\top$-NARE associated with the eigenvalues in $\mathcal U$.

First, we state a simple result that will be useful in the following.

\begin{lemma}\label{lem:u1inv}
Let $\varphi(z)=\mathcal M+z\mathcal M^\top$ be a regular matrix pencil of size $2n$ and let 
\[
    (\mathcal M+z\mathcal M^\top)\begin{bmatrix} U_1\\U_2\end{bmatrix} = \begin{bmatrix} V_1\\V_2\end{bmatrix}\alpha(z),
\]
where $U_1,U_2,V_1,V_2$ and $\alpha(z)$ are of size $n$, $\left[\begin{smallmatrix} U_1\\U_2\end{smallmatrix}\right]$ and $\left[\begin{smallmatrix} V_1 \\V_2\end{smallmatrix}\right]$ have rank $n$, and $\alpha(z)$ collects a reciprocal-free subset of eigenvalues of $\varphi(z)$. The matrix $U_1$ is invertible if and only if the matrix $V_2$ is invertible.
\end{lemma}
\begin{proof}
Let $U_1$ be invertible and $X:=U_2U_1^{-1}$. Since the eigenvalues of $\alpha(z)$ are a reciprocal-free set, the subspace spanned by the columns of $\left[\begin{smallmatrix} I\\X\end{smallmatrix}\right]$ is isotropic and
\[ 
\varphi(z)\begin{bmatrix} I\\X\end{bmatrix}
=\begin{bmatrix} V_1\\V_2\end{bmatrix} 
\alpha(z)U_1^{-1}.
\]
Premultiplying by $\left[\begin{smallmatrix} I & X^\top\end{smallmatrix}\right]$ the left hand side is $0$ and then
\[
	(V_1+X^\top V_2)\alpha(z)U_1^{-1}=0.
\]
Since $\alpha(z)$ is regular, we have that $V_1+X^\top V_2=0$. If $V_2$ is singular then $V_2v=0$ for $v\ne 0$ and thus $V_1v=0$ and this contradicts the fact that $\left[\begin{smallmatrix} V_1\\V_2\end{smallmatrix}\right]$ has dimension $n$. Thus $V_2$ is invertible.

The converse can be proved in a similar way by assuming that $V_2$ is invertible and writing
\[
	\varphi(z)\begin{bmatrix} U_1\\ U_2\end{bmatrix}
	=\begin{bmatrix} -X^\top\\I\end{bmatrix} V_2\alpha(z),
\]
with $X:=-(V_1V_2^{-1})^\top$. {}From isotropy we obtain $-U_1^\top X^\top+U_2^\top=0$, from which we deduce that $U_1$ is invertible.
\end{proof}

\begin{theorem}[Sufficient condition for invertibility of $U_1$\textbf{}]
Let a $\top$-NARE~\eqref{eq:tnare} be given, and let the columns of $\begin{bsmallmatrix}
U_1 \\ U_2\end{bsmallmatrix}$ span a deflating subspace associated with a reciprocal-free set $\mathcal{S}$ of eigenvalues of $\varphi(z)$.

Suppose that for each $v\in \mathbb{C}^n$ there exists $z_*\in\mathbb{C}$ such that
\[
	v^\top(B+z_* B^\top)v \neq 0.
\]
Then, $U_1$ is invertible, and the $\top$-NARE has a solution $X$ associated with the eigenvalues in $\mathcal{S}$. 
\end{theorem}
\begin{proof}
Suppose by contradiction that $U_1$ is singular; then there exists $v\ne 0$ such that $U_1v=0$. Set $v_2:=U_2v$; then $v_2 \neq 0$ because $\left[\begin{smallmatrix}U_1\\U_2\end{smallmatrix}\right]$ has full rank. From \eqref{eq:20} and isotropy we have that
\[
	v^\top\begin{bmatrix} U_1^\top & U_2^\top\end{bmatrix} \varphi(z_*)
	\begin{bmatrix} U_1\\U_2\end{bmatrix}v=0,
\]
that can be written as $v_2^\top(B+z_* B^\top)v_2=0$. This contradicts our assumption.
\end{proof}

The previous results implies that if $B$ is stable, that is $B+B^\top$ is definite positive, then a solution associated with a reciprocal-free subset of eigenvalues of $\varphi(z)$ exists.

\begin{theorem}
With the notation of the section, a necessary condition for $U_1$ to be invertible is that for any $z\not\in\mathcal{S}$ the matrix
\[
	\begin{bmatrix} A+zD^\top & -B-zB^\top\end{bmatrix}
\]
has rank $n$.
\end{theorem}
\begin{proof}
Let $u\ne 0$ be such that $u^\top (A+zD^\top)=u^\top(B+zB^\top)=0$ for some $z\not\in\mathcal{S}$.
We get
\[
	\begin{bmatrix} 0 & u^\top\end{bmatrix}\varphi(z)=0.
\]
Premultiplying equation \eqref{eq:20} by $\left[\begin{smallmatrix} 0 & u^\top\end{smallmatrix}\right]$ we get 
\[
	u^\top V_2\wt \alpha(z)=0,
\]
and since $z\not\in \mathcal S$ we can simplify $\wt \alpha(z)$ obtaining that $V_2$ is singular and in turn $U_1$ is singular in view of Lemma~\ref{lem:u1inv}.
\end{proof}

\section{Classical algorithms}\label{sec:cl}

In prior work \cite{bp20,bimp} some algorithms for the solution of $\top$-Riccati equations have been introduced. All of them, but Newton's method, rely on the linearization \eqref{eq:phial} and the solution of the resulting linear eigenvalue problem.

\paragraph{Newton's method.} 
Applying Newton's method to equation $\mathcal R(X)=0$ in \eqref{eq:tnare} one obtains the sequence
\[
    X_{k+1}=X_k+H_k,\qquad \mathcal {\mathcal R}'(X_k)[H_k]=-\mathcal R(X_k),
\]
for a given $X_0\in\mathbb R^{n\times n}$. The Newton step can be explicitly written as
\[
    X_{k+1}=X_k+H_k,\qquad (D-X_k^\top B)H_k+H_k^\top (A-BX_k)=-\mathcal R(X_k),
\]
and thus each step requires solving a $\top$-Sylvester equation. This results in an overall expensive algorithm, whose convergence to the required solution is hardly predictable, unless some very strict conditions are given on the coefficients \cite{bp20}.

\paragraph{QZ.} 
The deflating subspace associated with a given subset of reciprocal-free eigenvalues of the pencil $\varphi(z)=\mathcal{M}+z\mathcal{M}^\top \in\mathbb{R}^{2n\times 2n}$ can be obtained by putting the pair $(\mathcal{M},\mathcal{M}^\top )$ in upper (quasi-)triangular form with the QZ algorithm and then reordering the (block) diagonal so that the required eigenvalues lie in the left half diagonal entries \cite{gvl}.

In practice, a pair of unitary (orthogonal) matrices $Q$ and $Z$ is constructed such that $QMZ^*=A$ and $QM^\top Z^*=B$ are upper (quasi-)triangular and by blocking 
\begin{equation} \label{ABpencil}
Q\varphi(z)Z^* = 
    A+zB=\twotwo{A_{11}+zB_{11}}{A_{12}+zB_{12}}{0}{A_{22}+zB_{22}}
\end{equation}
using $n\times n$ blocks, the spectrum of $A_{11}+zB_{11}$ is the one required. Then the first $n$ columns of $Q$ spans the required deflating subspace of $\varphi(z)$.

The QZ algorithm \cite[Alg.~1]{bimp} spoils the $\top$-palindromic structure of the pencil, since the transformed pencil~\eqref{ABpencil} is not $\top$-palindromic.

\paragraph{Palindromic QZ.}
There exists a variant of the QZ algorithm that exploits the $\top$-palindromic structure of the pencil $\varphi(z)$ \cite{ksw09,m4}. The algorithm is based on a modification of the QZ algorithm which produces a $\top$-palindromic pencil $R+z R^\top $, where $R$ is an antitriangular matrix, that written with $n\times n$ blocks is
\[
    R=\twotwo{0}{R_{12}}{R_{21}}{R_{22}}.
\]

In practice, a unitary matrix $U$ is constructed such that $U^* \mathcal MU=R$ is antitriangular. The first $n$ columns of $U$ span a deflating subspace associated with the eigenvalues of $R_{21}+zR_{12}^\top$.

If the deflating subspace associated with a given reciprocal-free set of eigenvalues of $\varphi(z)$ is required, one must apply a transformation to $R+zR^\top $ that puts the required eigenvalues on the lower left block of the antitriangular pencil. Composing the two transformations one obtains $Q$
such that $Q^* \varphi(z)Q$ is antitriangular and the first $n$ columns of $Q$ span a deflating subspace associated with the required subset of the spectrum.
A procedure to move the required eigenvalues in the bottom left part of the antitriangular pencil is described in \cite[Sec. 5.3]{bimp}, but it requires $\mathcal{O}(n^4)$ arithmetic operations (ops) in the worst case. An $\mathcal{O}(n^3)$ procedure will be described in Section~\ref{sec:sw}.

\paragraph{Doubling.}
The doubling algorithm is an iterative procedure to get the deflating subspace of a pencil associated with eigenvalues inside the unit disk \cite[Alg.~2]{bimp}. It is applicable to the pencil $\varphi(z)$, when it is non-critical.

The algorithm is based on a transformation that takes as input a matrix pencil $\mathcal{A} + z\mathcal{B}$ and produces a new pencil $\widehat{\mathcal{A}} + z \widehat{\mathcal{B}}$ whose eigenvalues are the squares of those of the starting pencil, and with the same right deflating subspaces. Applying repeatedly this transformation to a non-critical pencil, one produces a sequence of pencils whose eigenvalues converge to either $0$ or $\infty$.

The same idea is then declined in different forms. The variant used in~\cite[Section~5.2]{bimp} uses a preprocessing step which produces a starting pencil right-equivalent to $\varphi(z)$ but in the so-called standard structured form
\[
    \twotwo{E_0}{0}{-P_0}{I}+z\twotwo{I}{-G_0}{0}{F_0}.
\]
The doubling iteration produces a sequence of pencils with the same structure, defined by parameters $\{(E_k,F_k,P_k,G_k)\}_{k}$. Then one has $\lim_k P_k=X$, where $X$ is the stabilizing solution.

\paragraph{Palindromic doubling.} A variant of the doubling algorithm that preserves the palindromic structure of the pencil is suggested in~\cite{pda}. The doubling transformation in this case takes as input a pencil of the form $\mathcal{M} + z \mathcal{M}^\top$, and produces the pencil $\widehat{\mathcal{M}} + z \widehat{\mathcal{M}}^\top$, with $\widehat{\mathcal{M}} = \mathcal{M}(\mathcal{M}+\mathcal{M}^\top)^{-1}\mathcal{M}$. The authors note that this transformation can be written in an alternative form: calling $\mathcal{H} = \frac{1}{2}(\mathcal{M} + \mathcal{M}^\top)$ and $\mathcal{K} = \frac{1}{2}(\mathcal{M} - \mathcal{M}^\top)$ the symmetric and anti-symmetric part of $\mathcal{M}$, respectively, the transformation becomes

\[
\widehat{\mathcal{H}} = \frac12 (\mathcal{H} + \mathcal{K}\mathcal{H}^{-1}\mathcal{K}), \quad \widehat{\mathcal{K}} = \mathcal{K}.
\]
In this form, one sees that this palindromic doubling algorithm is formally equivalent to applying a Cayley transform that transforms the initial pencil $\mathcal{M} + z\mathcal{M}^\top$ into $\mathcal{M}+\mathcal{M}^\top + z(\mathcal{M}-\mathcal{M}^\top) = 2\left(\mathcal{H} + z\mathcal{K}\right)$, and then applying to this pencil a classical iteration to compute the matrix sign function using symmetric matrices (see, e.g., Gardiner and Laub~\cite{gardiner}).

Since the algorithm uses the inverse of $\mathcal{H} = \frac12(\mathcal{M} + \mathcal{M}^\top)$, one can foresee numerical problems when the pencil $\varphi(z)$ has an eigenvalue close to $\lambda=1$.

\section{Swapping the eigenvalues in the antitriangular form of a $\top$-palindromic pencil}\label{sec:sw}

Once a $\top$-palindromic pencil $\varphi(z)=\mathcal{M}+z\mathcal{M}^\top$ has been reduced to an antitriangular form $\psi(z)=\mathcal{N}+z\mathcal{N}^\top$ by means of the palindromic QZ algorithm, the eigenvalues of the pencil can be read from the antidiagonal of $\psi(z)$. If a unitary matrix $U$ puts the pencil $\varphi(z)$ in antitriangular form (that is, $\mathcal{N}=U^* \mathcal{M}U$ is antitriangular), then the first $n$ columns of $U$ span an $n$-dimensional deflating subspace for the pencil $\varphi(z)$, associated with the eigenvalues appearing in the antidiagonal entries of the first $n$ columns.

In the applications, such as the solution of $\top$-NAREs, one is interested in a deflating subspace associated with a given subset of eigenvalues of $\varphi(z)$; for instance, the ones inside the unit disk. These eigenvalues do not necessarily appear in bottom left $n\times n$ block of the antidiagonal of $\psi(z)$; hence in order to select the required deflating subspace we need a procedure to reorder the antitriangular form, that is, a transformation putting the pencil $\psi(z)$ into a new antitriangular form, with the required eigenvalues in the bottom left block. A procedure for this purpose has been developed in \cite[Sec. 5.3]{bimp}, but unfortunately it requires $\mathcal{O}(n^4)$ arithmetic operations (ops) for a $2n\times 2n$ matrix pencil, in the worst case. Here we propose a procedure based on the swap of the spectrum of adjacent blocks on the block antidiagonal of a $\top$-palindromic antitriangular pencil yielding an algorithm that requires $\mathcal{O}(n^3)$ arithmetic operations in the worst case.

We consider operations of two types to perform the swap on an even-size $\top$-palindromic antitriangular pencil: the first one, described in Lemma \ref{thm:lemswap1}, swaps the spectrum of two blocks adjacent to the center of the antidiagonal; while the second, described in Lemma \ref{thm:lemswap2}, swaps the spectrum of two adjacent blocks in the upper part of the antidiagonal and at the same time the spectrum of the two corresponding adjacent blocks in the lower part.

In the following, we use the fact that the product $ABC$ is antitriangular if $A$ is lower triangular, $B$ is antitriangular and $C$ is upper triangular.

\begin{lemma}\label{thm:lemswap1}
Let
\[
	R=\begin{bmatrix}
	0 & R_{12}\\
	R_{21} & R_{22}
	\end{bmatrix}\in \mathbb C^{2m\times 2m},
\]
with $R_{12},R_{21}\in\mathbb C^{m\times m}$ antitriangular matrices such that $R_{21}+zR_{12}^\top $ is regular with  reciprocal-free spectrum. Then the matrix equation $R_{21}Y+Y^\top R_{12}=-R_{22}$ has a unique solution $Y$, and if
$
\left[\begin{smallmatrix}
Y & I \\
I & 0\\
\end{smallmatrix}\right]=QU
$ 
is a QR factorization, then
\[
Q^\top RQ = \begin{bmatrix}
0 & \wt R_{21}\\ \wt R_{12}& \wt R_{22}
\end{bmatrix},
\]
where $\wt R_{12},\wt R_{21}\in\mathbb C^{m\times m}$ are antitriangular and $\wt R_{12}+z\wt R_{21}^\top $ has the same eigevalues as $R_{12}+zR_{21}^\top $, moreover $\wt R_{21}+z\wt R_{12}^\top $ has the same eigenvalues as $R_{21}+zR_{12}^\top $.
\end{lemma}

\begin{proof}
Since $R_{21}+zR_{12}^\top $ has reciprocal-free spectrum, the matrix equation $R_{21}Y+Y^\top R_{12}=-R_{22}$ has a unique solution \cite{bk}. Then the proof can be obtained observing that
\[
\begin{bmatrix} Y^\top  & I \\ I & 0\end{bmatrix}R\begin{bmatrix} Y & I\\I & 0\end{bmatrix}=\begin{bmatrix} 0 & R_{21}\\R_{12} & 0\end{bmatrix},
\]
and that $Q^\top (R+zR^\top )Q=U^{-\top}\left[\begin{smallmatrix} 0 & R_{21}+zR_{12}^\top \\R_{12}+zR_{21}^\top  & 0\end{smallmatrix}\right]U^{-1}$ has the required property. For instance, the block of $Q^\top (R+zR^\top )Q$ in position $(1,2)$ is $(U^{-\top})_{11}(R_{21}+zR_{12}^\top )(U^{-1})_{22}$.
\end{proof}

To provide the next swapping technique, we need a characterization of the uniqueness of solution of a system of Sylvester equations that can be found in~\cite[Thm. 4.2]{stewart72} (see also \cite[Thm. 5]{dipr19}).
\begin{lemma}\label{thm:sylv}
Let $A,A'\in\C^{m\times m}$, $B,B'\in\C^{n\times n}$ and $C,C'\in\C^{m\times n}$. The system of matrix equations
\begin{equation}\label{eq:sylv}
\left\{\begin{array}{l}
AX-YB=C,\\ A'X-YB'=C',
\end{array}\right.
\end{equation}
admits a unique solution if and only if the pencils $A-\la A'$ and $B-\la B'$ are regular and with disjoint spectra.
\end{lemma}
Note that in the statement of~\cite[Theorem 4.2]{stewart72} one requires only that the pencils $A-\la A'$ and $B-\la B'$ have disjoint spectra, but this different formulation is due just to the fact that in \cite{stewart72} the spectrum is defined also for singular pencils (and in this case it is defined as~$\mathbb{C}\cup\{\infty\}$).

We are ready to describe the tool to swap the spectra of two pairs of blocks on the antidiagonal.

\begin{lemma}\label{thm:lemswap2}
Let
\[
	R=\begin{bmatrix}
	0 & R_{12}\\
	R_{21} & R_{22}
	\end{bmatrix}\in \mathbb C^{m\times m},\qquad
	S=\begin{bmatrix}
	0 & S_{12}\\
	S_{21} & S_{22}
	\end{bmatrix}\in \mathbb C^{m\times m},
\]
with $R_{12},S_{12}\in\mathbb C^{\ell\times \ell}$ and $R_{21},S_{21}\in\mathbb C^{(m-\ell)\times (m-\ell)}$ antitriangular matrices such that the pencils $R_{12}-\la S_{12}$ and $R_{21}-\la S_{21}$ are regular and with disjoint spectra.
Let $(X,Y)$ be the unique solution of 
the system of Sylvester equations
\[
\left\{\begin{array}{l}
XR_{12}+R_{21}Y=-R_{22},\\
XS_{12}+S_{21}Y=-S_{22}.\end{array}\right.
\]
If
$
\left[\begin{smallmatrix}
X^\top  & I \\
I & 0\\
\end{smallmatrix}\right]=QU
$ 
and 
$
\left[\begin{smallmatrix}
Y & I \\
I & 0\\
\end{smallmatrix}\right]=PV$ 
are QR factorizations, then
\[
Q^\top RP = \begin{bmatrix}
0 & \wt R_{21}\\ \wt R_{12}& \wt R_{22}
\end{bmatrix},\qquad
Q^\top SP = \begin{bmatrix}
0 & \wt S_{21}\\ \wt S_{12}& \wt S_{22}
\end{bmatrix},\qquad
\]
where $\wt R_{21},\wt S_{21}\in\mathbb C^{(m-\ell)\times (m-\ell)}$ and $\wt R_{12},\wt S_{12}\in\mathbb C^{\ell\times\ell}$ are antitriangular. Moreover, the eigenvalues of the pencils $\wt R_{12}+z\wt S_{12}$ and $\wt R_{21}+z\wt S_{21}$ are the same as the ones of $R_{12}+zS_{12}$ and $R_{21}+zS_{21}$, respectively.
\end{lemma}
\begin{proof}
Let $(X,Y)$ be the solution of the system of Sylvester equations, which is unique by Lemma~\ref{thm:sylv}. We have
\[
\begin{bmatrix} X & I \\ I & 0\end{bmatrix}(R+zS)\begin{bmatrix} Y & I\\I & 0\end{bmatrix}=\begin{bmatrix} 0 & R_{21}+zS_{21}\\R_{12}+zS_{12} & 0\end{bmatrix},
\]
from which we get
\[
	Q^\top (R+zS)P=U^{-\top}
	\begin{bmatrix} 0 & R_{21}+zS_{21}\\R_{12}+zS_{12} & 0\end{bmatrix}
	V^{-1}=\begin{bmatrix} 0 & \wt R_{21}+z\wt S_{21}\\ \wt R_{12}+z\wt S_{12} & \wt R_{22}+z\wt S_{22}\end{bmatrix},
\]
which is antitriangular, since $U$ and $V$ are triangular. Finally, we have that $\wt R_{21}+z\wt S_{21}=(U^{-\top})_{11} (R_{21}+zS_{21})(V^{-1})_{22}$ and a similar equation for $\wt R_{12}+z\wt S_{12}$, from which it follows that $\wt R_{21}+z\wt S_{21}$ and $\wt R_{12}+z\wt S_{12}$ have the same spectrum as $R_{21}+z S_{21}$ and $R_{12}+z S_{12}$, respectively. 
\end{proof}

The two previous lemmas allow one to swap the spectrum of blocks on the diagonal of an antitriangular $\top$-palindromic pencil. The first one allows us to swap the spectrum of two blocks of the same size adjacent to the center of the matrix, while the second to swap the spectrum of two blocks on the upper half of the matrix and the spectrum of corresponding two blocks on the lower half.

Assume that we want to swap the eigenvalues of $R_{12}+zR_{21}^\top \in\mathbb R^{m\times m}$ and $R_{21}+zR_{12}^\top \in\mathbb R^{m\times m}$, adjacent to the center of the $2n\times 2n$ pencil 
\[
\varphi(z)=M+zM^\top ,\qquad M=\begin{bmatrix}
 & & & \ast\\
 & & R_{12} & \ast\\
 & R_{21} & R_{22} & \ast\\
 \ast & \ast & \ast & \ast
\end{bmatrix},
\]
then it is sufficient to find $Q$ as in Lemma \ref{thm:lemswap1} and apply the congruence with matrix
\[
	V:=\begin{bmatrix} I_{n-m} \\ & Q\\ & & I_{n-m}\end{bmatrix}.
\]
to the pencil $\varphi(z)$.
Since $R_{12}$ has size $m$ and $\mathcal{M}$ has size $2n$, then the congruence $V^\top \mathcal{M}V$ modifies only $2m$ rows and columns (whose first $n-m$ elements are zeros) of $\mathcal M$ and thus it requires $16m^2(m+n)$ ops.

On the other hand, if we want to swap the eigenvalues of $R_{12}+zS_{12}\in\mathbb R^{\ell\times \ell}$ and $R_{21}+zS_{21}\in\mathbb R^{(m-\ell)\times (m-\ell)}$ of the pencil
\[
	\varphi(z)=\mathcal{M}+z\mathcal{M}^\top ,\qquad \mathcal M=\begin{bmatrix}
	& & & & & & \ast\\
	& & & & & R_{12} & \ast\\
	& & & & R_{21} & R_{22} & \ast\\
	& & & \ast & \ast & \ast & \ast & \\
	& & S_{21}^\top  & \ast & \ast & \ast & \ast\\
	& S_{12}^\top  & S_{22}^\top  & \ast & \ast & \ast & \ast\\
	\ast & \ast & \ast & \ast & \ast & \ast & \ast
	\end{bmatrix}
\]
we need to swap at the same time the eigenvalues of the pencils $S_{21}^\top +zR_{21}^\top $ and $S_{12}^\top +zR_{12}^\top $ with the congruence
\[
	V:=\begin{bmatrix}
	I_q \\ & Q \\ & & I_{2(n-m-q)} \\ & & & P\\ & & & & I_q
	\end{bmatrix},
\]
where $Q$ and $P$ are as in Lemma \ref{thm:lemswap2}.

Since the size of $R_{12}$ is $\ell$ and the size of $R_{21}$ is $m-\ell$, then the congruence $V^\top \mathcal{M}V$ modifies only $2m$ rows and columns ($m$ of them have only $m+q$ nonzero elements, while the other $m$ have $2n-q$ nonzero elements) in $M$ and thus it requires $4m^2(m+2n)$ ops.

If we are interested in the stable deflating subspace of a non-critical $\top$-palindromic pencil, then it is possible to move the eigenvalues with modulus smaller than one to the left half by using swaps involving only $1\times 1$ blocks on the diagonal.

Whenever an eigenvalue with modulus smaller than $1$ follows an eigenvalue with modulus greater than $1$, a swap is needed. If the two eigenvalues are adjacent to the center, then the procedure of Lemma \ref{thm:lemswap1} can be applied, while if the two eigenvalues are in one half of the matrix, then the procedure of Lemma \ref{thm:lemswap2} will swap the two eigenvalues and two more eigenvalues in the opposite part. Note that the hypotheses of Lemma \ref{thm:lemswap1} and \ref{thm:lemswap2} are fulfilled since the eigenvalues to swap are inside and outside the unit circle, respectively.

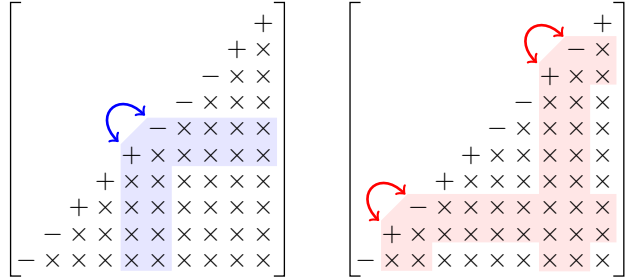
\begin{figure}
\centering
    \begin{tikzpicture}[scale=0.7]
    \draw[line width=0.5] (-0.60,-0.30) -- (-0.80,-0.30) -- (-0.80,4.90) -- (-0.60,4.90);
    \draw[line width=0.5] (4.20,-0.30) -- (4.40,-0.30) -- (4.40,4.90) -- (4.20,4.90);
    \filldraw[blue!10] (1.30,-0.20) -- (1.30,2.20) -- (1.80,2.70) -- (4.25,2.70) -- (4.25,1.80) -- (2.25,1.80) -- (2.25,-0.20);
    \foreach \i in {0,...,8}
        \foreach \j in {0,...,\i}
        {
            \node at (\i*0.5,\j*0.5) {$\times$};
        }
    \node at (-0.5,0) {$-$};
    \node at (0,0.5) {$-$};
    \node at (0.5,1) {$+$};
    \node at (1,1.5) {$+$};
    \node at (1.5,2) {$+$};
    \node at (2,2.5) {$-$};
    \node at (2.5,3) {$-$};
    \node at (3,3.5) {$-$};
    \node at (3.5,4) {$+$};
    \node at (4,4.5) {$+$};
    \draw[<->,color=blue,line width=1]
    (1.25,2.25) .. controls (0.75,2.75) and (1.25,3.25) .. (1.75,2.75);
    \end{tikzpicture}
    \qquad
            \begin{tikzpicture}[scale=0.7]
    \draw[line width=0.5] (-0.60,-0.30) -- (-0.80,-0.30) -- (-0.80,4.90) -- (-0.60,4.90);
    \draw[line width=0.5] (4.20,-0.30) -- (4.40,-0.30) -- (4.40,4.90) -- (4.20,4.90);
    \filldraw[color=red!10] (-0.2,-0.20) -- (-0.2,0.75) -- (0.30,1.25) -- (2.80,1.25) -- (2.80,3.75) -- (3.30,4.25) -- (4.25,4.25) -- (4.25,3.35) -- (3.75,3.35) -- (3.75,1.25) -- (4.25,1.25) -- (4.25,0.35) -- (3.75,0.35) -- (3.75,-0.20) -- (2.80,-0.20) -- (2.80,0.35) -- (0.75,0.35) -- (0.75,-0.20); 
    \foreach \i in {0,...,8}
        \foreach \j in {0,...,\i}
        {
            \node at (\i*0.5,\j*0.5) {$\times$};
        }
    \node at (-0.5,0) {$-$};
    \node at (0,0.5) {$+$};
    \node at (0.5,1) {$-$};
    \node at (1,1.5) {$+$};
    \node at (1.5,2) {$+$};
    \node at (2,2.5) {$-$};
    \node at (2.5,3) {$-$};
    \node at (3,3.5) {$+$};
    \node at (3.5,4) {$-$};
    \node at (4,4.5) {$+$};
    \draw[<->,color=red,line width=1] (-0.25,0.75) .. controls (-0.75,1.25) and (-0.25,1.75) .. (0.25,1.25);
    \draw[<->,color=red,line width=1]
    (2.75,3.75) .. controls (2.25,4.25) and (2.75,4.75) .. (3.25,4.25);
    \end{tikzpicture}
    \caption{Illustration of the swaps on an antitriangular matrix with highlighted the elements that are modified. Left: single swap; right: double swap.} \label{fig:swap}
\end{figure}

The swap procedure is depicted in Figure \ref{fig:swap} and described in Algorithm \ref{alg:swap}. A single swap consists in multiplying a $2\times 2$ matrix by the two central columns and rows of an antitriangular matrix, for a total cost of $12(n+1)$ ops. A double swap consists in multiplying $2\times 2$ matrices by two pairs of symmetric rows and columns, for a total cost of $24(n+1)$ ops.

In the worst case, that is, when all the required eigenvalues lie in the top right corner, $n$ single swaps and $\frac{1}{2}n(n-1)$ double swaps are required, hence the total cost of the procedure is $12n^2(n+1)$ ops (for a $2n\times 2n$ matrix).

\begin{algorithm}
\label{alg:swap}
\textbf{Input:} A unitary matrix $U$ and an antitriangular matrix $\mathcal{N}$ of size $2n$, such that the pencil $\mathcal{N}+z\mathcal{N}^\top$ is non-critical.

\textbf{Output:} A unitary matrix $Q$ and an antitriangular matrix $R$ such that $QRQ^\top=U\mathcal{N}U^\top$ and the first $n$ eigenvalues appearing on the antidiagonal of $R+zR^\top$ have modulus less than one.

Set $R\gets \mathcal{N}$, $Q\gets U$\;
\While{\rm
the set 
\[
    \Bigl\{1\le i\le n\,:\,
    \Bigl|\frac{R_{i,2n-i+1}}{R_{2n-i+1,i}}\Bigr|>
    1
    >\Bigl|\frac{R_{i+1,2n-i}}{R_{2n-i,i+1}}\Bigr|
    \Bigl\}
\]
is nonempty\;
}
{
Let $m$ be the minimum element of the set\;

\If{$m<n$}
{
 Solve
\[
    \twotwo{R_{m,2n-m+1}}{R_{m+1,2n-1}}{R_{2n-m+1,m}}{R_{2m-m,m+1}}
    \twoone{x}{y}
    =-\twoone{R_{m+1,2n-m+1}}{R_{2n-m+1,m+1}};\;
\]
Find the $Q$ factors $Q_x,Q_y$ of the QR factorizations of
$\begin{bsmallmatrix}{x} &{1} \\{1} & {0}\end{bsmallmatrix}$ and $\begin{bsmallmatrix}{y} &{1} \\{1} & {0}\end{bsmallmatrix}$, respectively\;
Update $R\gets V^\top RV$, $Q\gets QV$, with
\[
    V=\begin{bmatrix}
    I_{n-m}\\
    & Q_x\\
    & & I_{2m-4}\\
    & & & Q_y\\
    & & & & I_{n-m}
    \end{bmatrix},
\]
without forming $V$, but multiplying explicitly only the blocks $Q_x$ and $Q_y$\;
}
\If{$m=n$}
{
Set $y=-R_{m+1,m+1}/(R_{m,m+1}+R_{m+1,m})$\;
Find the $Q$ factor $Q_y$ of the QR factorization of $\begin{bsmallmatrix}{y} &{1} \\{1} & {0}\end{bsmallmatrix}$\;
Update $R\gets V^\top RV$, $Q\gets QV$ with
\[
    V=\begin{bmatrix}
    I_{n-1} \\ & Q_y\\ & & I_{n-1}
    \end{bmatrix},
\]
without forming $V$, but multiplying explicitly only the block~$Q_y$\;
}
}
\caption{Swap of the eigenvalues in the antitriangular form of a $\top$-palindromic pencil.}
\end{algorithm}

For completeness, the whole palindromic QZ algorithm for computing a stabilizing solution of the $\top$-NARE \eqref{eq:tnare} is summarized in  Algorithm~\ref{alg:palqz}.

\begin{algorithm}
\textbf{Input}: The matrix $\mathcal M$ of \eqref{eq:M}, such that the pencil $\varphi(z)=\mathcal M+z\mathcal{M}^\top$ is non-critical.

\textbf{Output}: An approximation $X$ to the stabilizing solution of  \eqref{eq:tnare}.

Compute a unitary matrix $U$ such that  $\mathcal N = U^\top \mathcal{M} U$ is antitriangular\;

Apply Algorithm~\ref{alg:swap} to $U$ and $\mathcal N$ and obtain $Q$\;

Return $X=Q_2 Q_1^{-1}$, where  $\begin{bsmallmatrix} Q_1 \\ Q_2 \end{bsmallmatrix}=Q \begin{bsmallmatrix}I_n\\ 0\end{bsmallmatrix} $.  

\caption{Palindromic QZ algorithm, with swapping, for computing the stabilizing solution the $\top$-NARE \eqref{eq:tnare}} \label{alg:palqz}
\end{algorithm}

\section{Quadraticization}\label{sec:qua}
In this section we show how the solutions of a $\top$-NARE \eqref{eq:tnare} can be related to the solutions of a suitable quadratic matrix equation.

We define the following quadratic matrix polynomial, obtained from $\varphi(z)$ by adding $n$ zeros at 0 and $n$ zeros at infinity:  
\begin{equation}\label{eq:Q}
	\mathcal Q(z):=
	\varphi(z)
	\begin{bmatrix} 0 & I \\ zI & 0\end{bmatrix}
	=\begin{bmatrix}
	Dz+A^\top z^2 & C+C^\top z\\
	-Bz-B^\top z^2 & A+D^\top z
	\end{bmatrix},
\end{equation}
and set $\mathcal Q(z)=\mathcal A_{-1}+\mathcal A_1 z+ \mathcal A_2 z^2$, so that
\begin{equation}\label{eq:at}
\mathcal A_{-1}=\begin{bmatrix} 0 & C \\ 0 & A\end{bmatrix},~~
\mathcal A_{0}=\begin{bmatrix} D & C^\top \\-B & D^\top \end{bmatrix},~~
\mathcal A_1=\begin{bmatrix} A^\top  & 0\\ -B^\top  & 0\end{bmatrix}.
\end{equation}
We may interpret $\mathcal Q(z)$ as a quadraticization of $\varphi(z)$.
A cubicization of $\varphi(z)$ is for instance the cubic matrix polynomial
\[
	\mathcal C(z):=\begin{bmatrix} I & 0\\0 & zI\end{bmatrix}
\varphi(z)\begin{bmatrix} zI & 0\\ 0 & I\end{bmatrix}
=\begin{bmatrix} -Bz-B^\top z^2 & A+D^\top z\\
Dz^2+A^\top z^3 & Cz+C^\top z^2\end{bmatrix}.
\]

The matrix coefficients of $\mathcal Q(z)$ and $\mathcal C(z)$ can be seen as block entries of suitable 
partitionings of the 
infinite dimensional block Toeplitz matrix
\begin{equation}\label{eq:bid}
\begin{bmatrix}
0 & \mathcal M & \mathcal M^\top  &  0 & 0 & \ldots \\
0 & 0 & \mathcal M & \mathcal M^\top  &  0 &  \ddots \\ 
\vdots& \ddots & \ddots & \ddots  & \ddots & \ddots 
\end{bmatrix},
\end{equation}
as shown in Figure~\ref{fig:part}.

\begin{figure}
\begin{center}
\begin{tikzpicture}[scale=0.8]
\draw node at (-2,0) {$0$} node at (-1,0) {$C$} node at (0,0) {$D$} node at (1,0.035) {$C^\top$} node at (2,0.035) {$A^\top $} node at (3,0) {$0$} node at (4,0) {$0$} node at (5,0) {$0$};
\draw node at (-2,-1) {$0$} node at (-1,-1) {$A$} node at (0,-1) {$-B$} node at (1,-1+0.035) {$D^\top$} node at (2,-1+0.035) {$-B^\top$} node at (3,-1) {$0$} node at (4,-1) {$0$} node at (5,-1) {$0$};
\draw node at (-2,-2) {$0$} node at (-1,-2) {$0$} node at (0,-2) {$0$} node at (1,-2) {$C$} node at (2,-2) {$D$} node at (3,-1.965) {$C^\top$} node at (4,-1.965) {$A^\top$} node at (5,-2) {$0$};
\draw node at (-2,-3) {$0$} node at (-1,-3) {$0$} node at (0,-3) {$0$} node at (1,-3) {$A$} node at (2,-3) {$-B$} node at (3,-2.965) {$D^\top$} node at (4,-2.965) {$-B^\top$} node at (5,-3) {$0$};
\draw[color=green!80!black,line width=0.8] (-1.35,0.35) -- (0.45,0.35) -- (0.45,-1.35) -- (-1.35,-1.35) -- (-1.35,0.35);
\draw[color=green!80!black,line width=0.8] (0.65,0.35) -- (2.5,0.35) -- (2.5,-1.35) -- (0.65,-1.35) -- (0.65,0.35);
\draw[color=blue!50,line width=0.8] (-0.3,0.3) -- (1.3,0.3) -- (1.3,-1.3) -- (-0.3,-1.3) -- (-0.3,0.3);
\draw[color=blue!50,line width=0.8] (-2.3,0.3) -- (-0.7,0.3) -- (-0.7,-1.3) -- (-2.3,-1.3) -- (-2.3,0.3);
\draw[color=blue!50,line width=0.8] (1.6,0.3) -- (3.3,0.3) -- (3.3,-1.3) -- (1.6,-1.3) -- (1.6,0.3);
\draw[color=orange,line width=0.8] (-2.4,-0.6) -- (-0.6,-0.6) -- (-0.6,-2.4) -- (-2.4,-2.4) -- (-2.4,-0.6);
\draw[color=orange,line width=0.8] (-0.4,-0.6) -- (1.4,-0.6) -- (1.4,-2.4) -- (-0.4,-2.4) -- (-0.4,-0.6);
\draw[color=orange,line width=0.8] (1.55,-0.6) -- (3.4,-0.6) -- (3.4,-2.4) -- (1.55,-2.4) -- (1.55,-0.6);
\draw[color=orange,line width=0.8] (3.6,-0.6) -- (5.4,-0.6) -- (5.4,-2.4) -- (3.6,-2.4) -- (3.6,-0.6);
\end{tikzpicture}
\caption{The coefficients of $\varphi(z)$ (green), $\mathcal Q(z)$ (blue) and $\mathcal C(z)$ (orange) in the block Toeplitz matrix \eqref{eq:bid}.}\label{fig:part}
\end{center}
\end{figure}

We are in particular interested in quadraticizations, since in the literature there are efficient algorithms for solving quadratic matrix equations. Indeed, the next theorem relates the solutions of the $\top$-NARE with the solutions of the quadratic matrix equation defined by $\mathcal Q(z)$.

\begin{theorem}\label{thm:bea}
Assume that $\varphi(z)$ is non-critical. The eigenvalues of $\mathcal Q(z)$ in \eqref{eq:Q} are the eigenvalues of $\varphi(z)$, plus $n$ eigenvalues equal to 0, and $n$ eigenvalues at infinity.
Moreover, 
equation \eqref{eq:tnare} has a stabilizing solution $S$ if and only if
the quadratic matrix equation
\begin{equation}\label{eq:Q2}
\mathcal A_{-1}+
\mathcal A_0 Z+
\mathcal A_1 Z^2=0
\end{equation}
has a solution $G$ with spectral radius less than one. In addition,
\begin{equation}\label{eq:GS}
G=\begin{bmatrix} 0 & S \\ 0 & W\end{bmatrix}
\end{equation}
with $W=-(D^\top -B^\top S)^{-1} (A-B S)$.
\end{theorem}

\begin{proof}
Since $\det \mathcal Q(z)=(-1)^n z^n \det \varphi(z)$, the eigenvalues of $\mathcal Q(z)$ in \eqref{eq:Q} are the eigenvalues of $\varphi(z)$, plus $n$ eigenvalues equal to 0, and $n$ eigenvalues at infinity. In particular, with the convention that $1/0=\infty$ and $1/\infty=0$, the eigenvalues of $\mathcal Q(z)$ come in pairs $(\lambda,1/\lambda)$.

Assume that $S$ is a stabilizing solution of \eqref{eq:tnare}. Then, from \eqref{eq:phial}, we have that
 $\det(D^\top -B^\top S)\ne 0$ since $\alpha(z)$ has not eigenvalues at infinity. According to \cite[Corollary 1]{bimp}, the matrix $S$ satisfies the equation
\begin{equation}\label{eq:dre}
C+DS+(C^\top +A^\top S)W=0,~~W=-(D^\top -B^\top S)^{-1} (A-B S).
\end{equation}
On the other hand,
by replacing the matrix $G$ of \eqref{eq:GS} in equation \eqref{eq:Q2}, we find that $G$ solves \eqref{eq:Q2} if and only if
\[
C+DS+C^\top  W + A^\top  S W=0, ~~ A-BS+D^\top  W- B^\top  S W=0.
\]
Such equations are satisfied since \eqref{eq:dre} holds true. Moreover, $G$ is the solution with minimal spectral radius since $\rho(G)=\rho(W)<1$, where the latter inequality follows from the definition of $\alpha(z)$ in \eqref{eq:phial}.

Conversely, let $G$ be a solution to \eqref{eq:Q2} such that $\rho(G)<1$. 
In particular $G$ has (at least) $n$ eigenvalues equal to 0. From the property $\mathcal A_{-1}\begin{bmatrix}
I \\ 0
\end{bmatrix}=\begin{bmatrix}
0 \\ 0
\end{bmatrix}$, it follows that $G$ is a matrix of the form 
\[
G=\begin{bmatrix}
I & \star \\
0 & \star 
\end{bmatrix}
\begin{bmatrix}
0 & \star\\
0 & \star
\end{bmatrix}
\begin{bmatrix}
I & \star \\
0 & \star 
\end{bmatrix}^{-1}=
\begin{bmatrix}
0 & G_1 \\
0 & G_2
\end{bmatrix},
\]
for suitable matrices $G_1$ and $G_2$, and where $\star$ stands for a generic $n\times n$ matrix. Since $G$ solves \eqref{eq:Q2}, we may easily verify that
\[
\mathcal Q(z)=(\mathcal A_0+\mathcal A_1 G+ \mathcal A_1 z)(zI-G).
\]
Therefore, since $\varphi(z)=\mathcal Q(z)\begin{bmatrix}
0 & z^{-1}I \\ I & 0
\end{bmatrix}$ and
\[
(zI-G)\begin{bmatrix}
0 & z^{-1}I \\ I & 0
\end{bmatrix}=
\begin{bmatrix}
zI & -G_1\\
0 & zI-G_2
\end{bmatrix}
\begin{bmatrix}
0 & z^{-1}I \\ I & 0
\end{bmatrix}
=
\begin{bmatrix}
-G_1 & I \\
zI-G_2 & 0
\end{bmatrix},
\]
we obtain
\[
\mathcal \varphi(z) \begin{bmatrix}
I \\ G_1
\end{bmatrix}=
(\mathcal A_0+\mathcal A_1 G+ \mathcal A_1 z)\begin{bmatrix}
0 \\ I
\end{bmatrix}
(zI-G_2).
\]
On the other hand, $\mathcal A_1 \begin{bmatrix}
0 \\ I
\end{bmatrix} =0$, therefore we conclude that
\[
\mathcal \varphi(z) \begin{bmatrix}
I \\ G_1
\end{bmatrix}=
\begin{bmatrix}
V_1 \\ V_2
\end{bmatrix} (zI-G_2),~~
\begin{bmatrix}
V_1 \\ V_2
\end{bmatrix}= 
(\mathcal A_0+\mathcal A_1 G)\begin{bmatrix}
0 \\ I
\end{bmatrix}.
\]
Since $\rho(G_2)=\rho(G)<1$, then the eigenvalues of $zI-G_2$ are a reciprocal-free set. Therefore, in view of \cite[Theorem 2]{bimp}, we conclude that $G_1$ solves \eqref{eq:tnare}. Moreover, according Lemma~\ref{lem:u1inv}, $V_2=D^\top -B^\top  G_1$ is invertible and, from \eqref{eq:phial}, we deduce that $G_2=W$.
\end{proof}

We can prove a similar result concerning the dual $\top$-NARE
\begin{equation}\label{eq:tnared}
  \mathcal S(Y)=0,\qquad \mathcal S(Y):=AY+ Y^\top D+Y^\top CY-B.
\end{equation}
In  particular, equation \eqref{eq:tnared} has a stabilizing solution $T$  if and only if
the quadratic matrix equation
$
\mathcal A_{-1} Z^2+
\mathcal A_0 Z+
\mathcal A_1 =0
$
has a solution $F$ with spectral radius less than one. Moreover, $T$ is the $(2,1)$ block of $F$.

According to Theorem~\ref{thm:bea}, if $\varphi(z)$ is non-critical, then the stabilizing solution $S$ of the $\top$-NARE \eqref{eq:tnare} can be computed from the solution  of \eqref{eq:Q2} with minimal spectral radius. To this purpose, a valid algorithm is the Cyclic Reduction algorithm \cite{bm:cr}. The overall procedure is summarized in Algorithm~\ref{algo:quad}. 
Concerning the computational cost, from the zero pattern of the matrices in \eqref{eq:at}, we deduce that the first block column of $\mathcal A_{-1}^{(k)}$ and the second block column of $\mathcal A_1^{(k)}$ are equal to zero, for $k=1,2,\ldots$. This property leads to a saving of the number of operations at each step, which is $\mathcal{O}(n^3)$ arithmetic operations. Concerning convergence properties, if Cyclic Reduction can be applied without breakdown and if the sequence $(\mathcal{\widehat A}^{(k)})^{-1}$ is uniformly bounded, then, in view of \cite[Theorems 3 and 4]{bm:cr}, we have
\[
\begin{split}
&\limsup_{k\to\infty} \sqrt[2^k]{ \Vert \mathcal A_{-1}^{(k)} \Vert}\le \sigma,~~
\limsup_{k\to\infty} \sqrt[2^k]{ \Vert \mathcal A_{1}^{(k)} \Vert }\le \sigma\\
&\limsup_{k\to\infty} \sqrt[2^k]{ \Vert \mathcal A_{0}^{(k)} - \mathcal H_0^{-1}\Vert }\le \sigma^2,~~
\limsup_{k\to\infty} \sqrt[2^k]{ \Vert 
G +(\widehat{\mathcal A}^{(k)})^{-1}\mathcal A_{-1} \Vert }
\le \sigma^2, \\
\end{split}
\]
where $\sigma=\max\{|z|\,:\,\det\varphi(z)=0,\, |z|<1\}$ and $\mathcal{H}_0$ is the constant coefficient of the Laurent matrix power series $\mathcal{H}(z)=(z^{-1}\mathcal{Q}(z))^{-1}$, $z\in\mathbb{T}$.

\begin{algorithm}
\textbf{Input}: matrices $\mathcal{A}_i$, $=-1,0,1$, maximum number of iterations $N$, tolerance $\varepsilon>0$ for the stopping condition.

\textbf{Output}: an approximation $X$ to the stabilizing solution of  \eqref{eq:tnare}.

Set $\mathcal A_i^{(0)}=\mathcal A_i$, $i=-1,0,1$, 
and $\widehat{\mathcal A}^{(0)}=\mathcal A_0$\;

\For{$k=1,2,\dots,N$}{
    Compute the matrices
    \[
    \begin{split}
        & \mathcal A_i^{(k+1)}=-\mathcal A_i^{(k)}(\mathcal A_0^{(k)})^{-1}\mathcal A_i^{(k)}\quad \text{for }i=-1,1\\
        & \mathcal A_0^{(k+1)}=\mathcal A_0^{(k)}-\mathcal A_{-1}^{(k)}(\mathcal A_0^{(k)})^{-1}\mathcal A_1^{(k)}- \mathcal A_{1}^{(k)}(\mathcal A_0^{(k)})^{-1}\mathcal A_{-1}^{(k)}\\
        & \widehat{\mathcal A}^{(k+1)}=\widehat {\mathcal A}^{(k)}-\mathcal A_{1}^{(k)}(\mathcal A_0^{(k)})^{-1}\mathcal A_{-1}^{(k)}\;
    \end{split}
    \]
    \If{$\Vert \mathcal A_1^{(k+1)}\Vert_\infty<\varepsilon \Vert \mathcal A_0^{(k+1)}\Vert_\infty$}{exit\;}
}
$G \gets -(\widehat{\mathcal A}_0^{(k+1)})^{-1}\mathcal A_{-1}$\;
Return $X$ from the $(1,2)$ block of $G$.

\caption{Cyclic Reduction for computing the stabilizing solution of the $\top$-NARE \eqref{eq:tnare}} \label{algo:quad}
\end{algorithm}

Another quadraticization can be obtained by 
defining the quadratic matrix polynomial $\mathcal S(z)=\varphi(-z)^2$, so that
\begin{equation}
\mathcal S(z)=\mathcal M^2-(\mathcal M \mathcal M^\top  + \mathcal M^\top  \mathcal M)z+
(\mathcal M^\top )^2 z^2.
\end{equation}
Observe that, unlike $\mathcal Q(z)$, $\mathcal S(z)$ is a $\top$-palindromic matrix polynomial. Moreover, the roots of $\mathcal S(z)$ are equal to the roots of $\mathcal \varphi(-z)$, with doubled multiplicity.

Theorem 5 of \cite{bimp} shows that, if $\varphi(z)$ is non-critical and has a stable deflating subspace with a basis of the form $\begin{bmatrix}
I \\ X
\end{bmatrix}$, then
\begin{equation}\label{eq:is}
\mathcal M\begin{bmatrix}
I \\ X
\end{bmatrix}
   = \mathcal M^\top
   \begin{bmatrix}
   I \\ X
   \end{bmatrix} W
\end{equation}
where $W=(D^\top - B^\top X )^{-1} (A - B X )$ and $\rho(W)<1$.

From \eqref{eq:is}, we deduce that
\[
\mathcal M^2 \begin{bmatrix}
I \\ X
\end{bmatrix}
   = 
  \mathcal M \mathcal M^\top 
   \begin{bmatrix}
   I \\ X
   \end{bmatrix} W, ~~ \mathcal M^\top  \mathcal M \begin{bmatrix}
I \\ X
\end{bmatrix}
   = 
  (\mathcal M^\top )^2
   \begin{bmatrix}
   I \\ X
   \end{bmatrix} W.
\]
Therefore, we obtain
\[
\mathcal M^2 \begin{bmatrix}
   I \\ X
   \end{bmatrix}
-(\mathcal M \mathcal M^\top  + \mathcal M^\top  \mathcal M) \begin{bmatrix}
   I \\ X
   \end{bmatrix} W +
(\mathcal M^\top )^2 \begin{bmatrix}
   I \\ X
   \end{bmatrix} W^2=0.
\]
By setting 
\[
\mathcal C=\mathcal M^\top  \mathcal M,~~
\mathcal B_{-1}=-\mathcal M^2, ~~ \mathcal B_0= \mathcal M^\top  \mathcal M +
\mathcal M \mathcal M^\top ,~~\widehat X=\begin{bmatrix}
   I \\ X
   \end{bmatrix},
\]
the above equalities lead to the homogeneous infinite dimensional system
\begin{equation}\label{eq:infs}
\begin{bmatrix}
\mathcal C & \mathcal{B}_{-1}^{\top}\\
\mathcal{B}_{-1} & \mathcal B_0 & \mathcal{B}_{-1}^{\top}\\
 &   \mathcal{B}_{-1} & \mathcal B_0 & \ddots\\
 && \ddots & \ddots 
\end{bmatrix}
\begin{bmatrix}
\widehat X \\
\widehat X W \\
\widehat X W^2 \\
\vdots
\end{bmatrix}=0.
\end{equation}
Cyclic Reduction applied to the above system generates the sequences of matrices
\begin{equation}\label{eq:cr}
    \begin{split}
   & \mathcal{B}_{-1}^{(k+1)}=  - \mathcal{B}_{-1}^{(k)} ({\mathcal{B}_{0}}^{(k)})^{-1} \mathcal{B}_{-1}^{(k)}  \\
   & \mathcal V^{(k+1)}= ({\mathcal{B}_{-1}}^{(k)})^\top  ({\mathcal{B}_{0}}^{(k)})^{-1} {\mathcal{B}_{-1}}^{(k)} \\
    & {\mathcal{C}}^{(k+1)}=   {\mathcal{C}}^{(k)} - {\mathcal{V}}^{(k+1)}\\
   & {\mathcal B_0}^{(k+1)}={\mathcal B_0}^{(k)}- {\mathcal{V}}^{(k+1)}- ({\mathcal{V}}^{(k+1)})^\top ,~~k=0,1,\ldots\\
    \end{split}
\end{equation}
with $\mathcal B_{-1}^{(0)}=\mathcal B_{-1}$, $\mathcal B_{1}^{(0)}=\mathcal B_{1}$, $\mathcal C^{(0)}=\mathcal C$, 
such that
\[
\mathcal{B}_{-1}^{(k)} \widehat X
  +
\mathcal{B}_0^{(k)} \widehat X W^{2^k} +
(\mathcal{B}_{-1}^{(k)})^\top  \widehat X W^{2^{2k}}=0
\]
and
\[
\mathcal{C}^{(k)} \widehat X + (\mathcal{B}_{-1}^{(k)})^\top  \widehat{X} W^{2^{2k}}=0. 
\]

Under the assumption that $\varphi(z)$ is non-critical, if $\det\mathcal M\ne 0$, then $\det(z^{-1}\mathcal S(z))\ne 0$ for $z\in\mathbb T$. Therefore,
according to the results of \cite[Section 2.2]{bea:extr} and \cite{bm:cr}, the sequences $\{\mathcal B_{-1}^{(k)}\}_k$ and $\{\mathcal B_{1}^{(k)}\}_k$ are well defined and convergent, moreover
$\limsup_{k\to\infty} \sqrt[2^k]{ \Vert \mathcal B_{-1}^{(k)} \Vert}\le \sigma$ and
$\limsup_{k\to\infty} \sqrt[2^k]{ \Vert \mathcal B_{0}^{(k)} - \mathcal K_0^{-1}\Vert }\le \sigma^2$, where $\mathcal K_0$ is the constant coefficient of the Laurent matrix power series $(z^{-1}\mathcal S(z))^{-1}$, for $|z|=1$. Similarly, we may prove that
the sequence $\{\mathcal C^{(k)}\}_k$ has a limit $\mathcal C^*$. Therefore, since $\rho(W)<1$, we conclude that $\mathcal C^* \widehat X=0$. Hence,
if the kernel of $\mathcal C^*$ has dimension $n$, then $\widehat X$ can be recovered from the null space of $\mathcal C^*$. The resulting procedure is summarized in Algorithm~\ref{algo:crs}. The computational cost of Algorithm~\ref{algo:crs} is roughly half the computational cost of Algorithm~\ref{algo:quad}, thanks to the symmetry of the involved matrices. However, from the numerical experiments,  Algorithm~\ref{algo:crs} seems to be less robust, in terms of accuracy, than Algorithm~\ref{algo:quad}.

\begin{algorithm}
\textbf{Input}: matrices $\mathcal B_i$, $i=-1,0$, and $\mathcal C$, maximum number of iterations $N$, tolerance $\varepsilon>0$ for the stopping condition.

\textbf{Output}: an approximation $X$ to the stabilizing solution of  \eqref{eq:tnare}.

Set $\mathcal B_i^{(0)}=\mathcal B_i$, $i=-1,0$, 
and ${\mathcal C}^{(0)}=\mathcal C$\;

\For{$k=1,2,\dots,N$}{
    Apply formulas \eqref{eq:cr}\;
    \If{$\Vert \mathcal B_{-1}^{(k+1)}\Vert_\infty<\varepsilon \Vert \mathcal B_{0}^{(k+1)}\Vert_\infty$}{exit\;}
}
Compute a rank-revealing QR decomposition $QRP = \mathcal C^{(k+1)}$\;
Return $X = Q_2 Q_{1}^{-1}$, where $\begin{bsmallmatrix} Q_1 \\ Q_2 \end{bsmallmatrix}=Q \begin{bsmallmatrix}I_n\\ 0\end{bsmallmatrix} $.

\caption{Symmetric Cyclic Reduction for computing the stabilizing solution to the $\top$-NARE \eqref{eq:tnare}} \label{algo:crs}
\end{algorithm}

\section{Integral representation}\label{sec:ir}

Another numerical method that benefits from the symmetry structure of the problem can be obtained by considering complex contour integrals. Assume that $\varphi(z)$ is invertible on the unit circle, and consider the contour integral
\[
\mathcal{I}(\mathcal{M}) = \frac{1}{2\pi \iunit}\oint \varphi(z)^{-1} \, dz,
\]
where the integral is taken on the unit circle $\mathbb{T}$. The following result holds.
\begin{theorem}
Assume that $\varphi(z)$ of \eqref{eq:phi} is non-critical. Then, $\operatorname{Range} \mathcal{I}(\mathcal{M})$ is the stable deflating subspace of $\varphi(z)$. Moreover,
    \[
    \mathcal{I}(\mathcal{M}) \mathcal{M}^\top
    \]
is the orthogonal projector on the stable deflating subspace.
\end{theorem}
\begin{proof}
To evaluate the integral, let us use a decomposition into Kronecker canonical form~\cite{gant} $\varphi(z) = P J(z) Q^{-1}$. Since $\varphi(z)$ is regular, $J(z)$ is the direct sum of Jordan blocks of the form $J_\lambda(z) = (\lambda-z)I + S$ with possibly Jordan blocks at infinity of the form $J_\infty(z) = I - zS$ (with $S$ the upper shift matrix with $S_{i,i+1}=1$ for each $i$ and all other elements equal to $0$). We shall further assume that this canonical form is ordered so that the leading $n\times n$ block of $J(z)$ contains the blocks $J_\lambda$ with $\abs{\lambda}<1$.

Then,
\begin{equation} \label{int1}
\mathcal{I}(\mathcal{M}) = Q\left(\frac{1}{2\pi \iunit}\oint J(z)^{-1}\, dz\right) P^{-1}.    
\end{equation}
A standard computation using the Cauchy integral formula gives
\[
\frac{1}{2\pi \iunit}\oint J_\lambda(z)^{-1}\, dz = \begin{cases}
-I & \abs{\lambda} < 1,\\
0 & \abs{\lambda} > 1,
\end{cases}
\quad
\frac{1}{2\pi \iunit}\oint J_\infty(z)^{-1}\, dz = 0.
\]
Hence,
\begin{equation}
\frac{1}{2\pi \iunit}\oint J(z)^{-1}\, dz = \begin{bmatrix}-I_n & 0\\ 0 & 0\end{bmatrix}.
\end{equation}
Thus, 
\begin{equation} \label{Imformula}
    \mathcal{I}(\mathcal{M}) = Q\begin{bmatrix}
    -I_n & 0\\
    0 & 0
    \end{bmatrix}P^{-1},
\end{equation}
and 
$\operatorname{Range} \mathcal{I}(\mathcal{M}) = \operatorname{Range} Q \begin{bsmallmatrix}I_n\\0\end{bsmallmatrix}$ is the span of the Jordan chains relative to the stable eigenvalues, which is precisely the stable deflating subspace.

By considering the degree-1 part of $\varphi(z) = P J(z) Q^{-1}$, one sees that
\begin{equation} \label{int2}
\mathcal{M}^\top = P \begin{bmatrix}-I_n & 0\\ 0 & \star\end{bmatrix} Q^{-1}
\end{equation}
(with $\star$ a block whose content is not relevant here). Multiplying~\eqref{Imformula} and~\eqref{int2}, we get
\[
 \mathcal{I}(\mathcal{M}) \mathcal{M}^\top = Q \begin{bmatrix}
I_n & 0\\
0 & 0
\end{bmatrix}Q^{-1},
\]
which is the second part of the thesis.
\end{proof}
Once the projector $\Pi:=\mathcal I(\mathcal M)\mathcal M^\top$ is computed, its range is the stable invariant subspace, so we can compute a basis $\begin{bsmallmatrix}U_1 \\ U_2\end{bsmallmatrix}$ for it and get the stabilizing solution as $X = U_2 U_1^{-1}$.

One can get approximations of $\mathcal I(\mathcal{M})$ with the trapezoidal rule on $2^k$ nodes
\[
\mathcal{I}(\mathcal{M}) \approx \mathcal{I}_{2^k}(\mathcal{M}) = \frac{1}{2^k}\sum_{j=1}^{2^k} \zeta^j \varphi(\zeta^j)^{-1},
\]
where $\zeta = \exp(\frac{2\pi \iunit}{2^k})$ is the principal $2^k$-th complex root of 1. Note that
\[
\zeta^j \varphi(\zeta^j)^{-1} = \zeta^{j/2} (\zeta^{-j/2} \mathcal{M} + \zeta^{j/2}\mathcal{M}^\top)^{-1}, 
\]
and the matrix in parentheses is Hermitian; hence the computation of $\mathcal{I}(\mathcal{M})$ can exploit the symmetry structure of the pencil. A full algorithm is presented in Algorithm~\ref{algo:trap}.
\begin{algorithm}
\textbf{Input}: $\mathcal{M} \in \mathbb{R}^{2n\times 2n}$ defined in \eqref{eq:M} such that $\varphi(z)=\mathcal{M}+z \mathcal{M}^\top$ is non-critical, number of nodes $2^k$.\\
\textbf{Output}: 
an approximation $X$ to the stabilizing solution of  \eqref{eq:tnare}.\\
$S \gets 0$\;
\For{$j=1,2,3,\dots,2^k$}{
    Compute the inverse of the Hermitian matrix $\zeta^{-j/2} \mathcal{M} + \zeta^{j/2}\mathcal{M}^\top$\;
    $S \gets S + \zeta^{j/2}(\zeta^{-j/2} \mathcal{M} + \zeta^{j/2}\mathcal{M}^\top)^{-1}$\;
}
$\Pi \gets \frac{1}{2^k} S \mathcal{M}^\top$\;
Compute a rank-revealing QR decomposition $QRP = \Pi$\;
$
\begin{bsmallmatrix}U_1 \\ U_2\end{bsmallmatrix} \gets Q \begin{bsmallmatrix}I_n\\ 0\end{bsmallmatrix}$\;
$X\gets U_2 U_1^{-1}$\;
\caption{Trapezoidal rule algorithm for computing the stabilizing solution of the $\top$-NARE \eqref{eq:tnare}} \label{algo:trap}
\end{algorithm}

If $\mathcal{M}$ is invertible, one can further simplify
\begin{equation} \label{closedintegral}
\mathcal{I}_{2^k}(\mathcal{M})\mathcal{M}^\top = \frac{1}{2^k}\sum_{j=1}^{2^k} \zeta^j(\mathcal{M}^{-\top}\mathcal{M} + \zeta^j I)^{-1} = (I-(-\mathcal{M}^{-\top}\mathcal{M})^{2^k})^{-1}.
\end{equation}
The last step follows from the scalar identity
\[
\frac{1}{2^k}\sum_{j=1}^{2^k} \frac{\zeta^j}{\zeta^j-x} = \frac{1}{1-x^{2^k}},
\]
as the same identity must hold for the corresponding matrix functions~\cite[Theorem~1.16]{highamfun}.

The formula~\eqref{closedintegral} reveals connections with doubling algorithms: comparing it with formulas such as~\cite[Equation~(4.8)]{BaiDG97} reveals that (in exact arithmetic) this method computes the same approximation of the stable invariant subspace as $k$ steps of the doubling algorithm. On the other hand, the cost of this algorithm is much larger than that of doubling algorithms, since it scales as $\mathcal{O}(2^k n^3)$ rather than $\mathcal{O}(kn^3)$.

Another alternative way in which Algorithm~\ref{algo:trap} can be interpreted is the evaluation-interpolation technique that appears, for instance, in \cite[Section~8.5]{blm}. Given a function $\psi(z)$ (scalar- or matrix-valued) that can be written as a Laurent power series $\psi(z) = \sum_{k\in\mathbb{Z}} c_k z^k$, converging in an annulus that contains the unit circle, we denote by $[\psi(z)]_k = c_k$ the coefficient of $z^k$ in it. Note that $\varphi(z)$ has no zeros on the unit circle, so its inverse $\psi(z) = \varphi(z)^{-1}$ can be written in this form. We have the following result.
\begin{lemma}
With the notation of the section, we have 
\[
[\varphi(z)^{-1}]_{-1} = \mathcal{I}(\mathcal{M}).
\]
\end{lemma}
\begin{proof}
We rely on a Kronecker form $\varphi(z) = PJ(z) Q^{-1}$ as above. Each block of $J(z)^{-1}$ has Laurent power series
\begin{equation} \label{Jordanlaurent}
J_\lambda(z)^{-1} = \begin{cases}
\sum_{j=0}^{\infty} (\lambda I + S)^{-j-1}z^j & \abs{\lambda} > 1,\\
\sum_{j=0}^{\infty} -(\lambda I + S)^{j}z^{-j-1} & \abs{\lambda} < 1.
\end{cases}
\end{equation}
Note, indeed, that the $j$\textsuperscript{th} coefficient $[J_\lambda(z)^{-1}]_j$ converges to zero when $j\to \pm\infty$, so these sums converge. For the blocks at infinity we obtain instead a finite sum
\[
{J_{\infty}(z)}^{-1} = (I-zS)^{-1} = \sum_{j=0}^{\infty} S^j z^j = \sum_{j=0}^{n-1} S^j z^j.
\]

By combining these series for each Jordan block, one obtains that the coefficient of $z^{-1}$ in the Laurent power series for $J(z)^{-1}$ is
\[
[J(z)^{-1}]_{-1} = \begin{bmatrix}
-I_n & 0\\
0 & 0
\end{bmatrix}.
\]
Comparing with~\eqref{Imformula} above, one gets the required result.
\end{proof}
We recall that, given a (possibly matrix-valued) function $\psi(z)$, its inverse discrete Fourier transform with $N$ nodes is the list of its evaluations at the roots of unity
\[
L = IDFT_{2^k}(\psi(z)) = [\psi(1), \psi(\zeta), \dots , \psi(\zeta^{2^k-1})],
\]
and that, given such a list of evaluations $L$, its discrete Fourier transform $DFT_{2^k}(L)$ is the unique polynomial  $p(z) = \sum_{j=0}^{2^k-1} c_j z^j$ of degree at most $2^k-1$ such that $IDFT_{2^k}(p(z)) = L$. The following result shows that we can obtain the coefficient of $z^{-1}$ in a Laurent power series as the limit of a certain evaluation/interpolation procedure using discrete Fourier transforms.
\begin{lemma}
Let $\psi(z)$ be a Laurent series converging on the unit circle, and $L = IDFT_{2^k}(\psi(z))$. Then,
\[
\lim_{k\to\infty} [DFT_{2^k}(L)]_{2^k-1} = [\psi(z)]_{-1},
\]
\end{lemma}
\begin{proof}
The result follows from~\cite[Theorem~3.8]{blm}.
\end{proof}
These two results suggest a method to compute an approximation of $\mathcal I(\mathcal{M})$: take a sufficiently large $k$, and compute
\[
L = IDFT_{2^k}(\varphi(z)^{-1}) = [\varphi(1)^{-1}, \varphi(\zeta)^{-1}, \dots , \varphi(\zeta^{2^k-1})^{-1}],
\]
directly by evaluating $\varphi(z)$ at the roots of unity and inverting the resulting matrices; then take the last coefficient of its DFT
\[
\mathcal{I}(\mathcal{M}) = [\varphi(z)^{-1}]_{-1}  \approx [DFT_{2^k}(L)]_{2^k-1}.
\]
In fact, this algorithm computes the same approximation $\mathcal{I}_{2^k}(\mathcal{M})$ obtained by the trapezoidal rule algorithm with $2^k$ nodes: indeed, by the matrix formulation of DFT, this last coefficient can be written as
\[
[DFT_{2^k}(L)]_{2^k-1} =  
\sum_{j=0}^{2^k-1} \overline{\zeta^{-j}} v_j =
\sum_{j=0}^{2^k-1} \overline{\zeta^{-j}} \varphi(\zeta^j)^{-1} = \mathcal{I}_{2^k}(\mathcal{M}).
\]

\section{Numerical experiments}\label{sec:exp}

The numerical experiments have the goal to compare the performances of the algorithms for solving the $\top$-NARE and to show the effectiveness of the swapping procedure proposed in Section~\ref{sec:sw}. The experiments have been performed on an Intel Core i5-1135G7 running Matlab R2021a on Ubuntu Linux 22.10.

\subsection{Numerical solution of the $\top$-NARE}
We test numerically  the various algorithms described in this paper to solve the $\top$-Riccati equation \eqref{eq:tnare}. In detail:
\begin{description}
\item[QZ] The (unstructured) QZ method;
\item[PalQZ] The palindromic QZ method (Algorithm~\ref{alg:palqz});
\item[DA] The (unstructured) doubling algorithm, as described in~\cite[Section~5.2]{bimp};
\item[CR1] Cyclic reduction with the non-symmetric quadraticization~\eqref{eq:Q} (Algorithm~\ref{algo:quad});
\item[CR2] Cyclic reduction on the symmetric infinite system~\eqref{eq:infs} (Algorithm~\ref{algo:crs});
\item[PDA] The palindromic doubling algorithm / sign function method from~\cite{pda}. We note that selecting the number of steps is harder with this method than with the other doubling-based methods, since the iterates are not guaranteed to converge up to the last significant digit in machine arithmetic --- unlike, for instance, Cyclic Reduction, in which $\norm{\mathcal{A}^{(k)}_{-1}}\norm{\mathcal{A}^{(k)}_1}$ converges to zero and hence at some point we obtain $\widehat{\mathcal{A}}^{(k+1)} = \widehat{\mathcal{A}}^{(k)}$ exactly. We have implemented a stopping criterion which obtains approximately the same number of steps as the other methods, and hence appears to be effective.
\item[Int] The integral representation method suggested in Section~\ref{sec:ir} (Algorithm~\ref{algo:trap}). For this method, again, selecting a stopping criterion is not trivial, as the iterates are not guaranteed to converge up to the last significant digit in machine arithmetic. For simplicity and ease of comparison, we just ran the same method with the same value of $k$ as \textbf{DA}.
\end{description}
For our experiments we use some of the test cases in~\cite{bimp}, and follow their numbering. The implementations of \textbf{QZ}, \textbf{PalQZ},
\textbf{DA} are the ones used in \cite{bimp}, apart from a minor variation in the stopping criterion for the last one.

\begin{description}
\item[Example 1] Example 1 in~\cite{bimp} with $n=10$, i.e.,
\begin{align*}
    A &= \begin{bmatrix}
    -1 & -1\\
    & -1 & -1\\\
    & & \ddots & \ddots\\
    & & & -1 & -1\\
    & & & & -1
    \end{bmatrix}, \quad
    D = 
    \begin{bmatrix}
    4 & -1\\
    & 4 & -1\\\
    & & \ddots & \ddots\\
    & & & 4 & -1\\
    & & & & 4
    \end{bmatrix},
    \\
    E &= \begin{bmatrix}
    -1 & -1\\
    & -1 & -1\\\
    & & \ddots & \ddots\\
    & & & -1 & -1\\
    & & & & -0.9
    \end{bmatrix},
\end{align*}
with $B = -A / \norm{A}_F$, $C = E / \norm{E}_F$. This is a small-scale, well-conditioned problem that should cause no difficulties to numerical algorithms.
\item[Example 2a] A larger example ($n=324$) in which $A,D$ are equal to the 2-dimensional finite difference stencil on a square, and $B,C$ are random matrices with coefficients in $[0,1]$ (obtained by the Matlaba command \verb!rand(n)!).
\item [Example 2b] The same example as above, but with $n=784$. This provides a larger-scale example, while still being solvable with a $\mathcal{O}(n^3)$ dense linear algebra algorithm.
\item[Example 4a] Example 4 in~\cite{bimp}, with $n=3$ and $\sigma = 10^{-10}$. This example has a reciprocal pair of eigenvalues of the form $\left( \frac{1}{1+\sigma}, 1+\sigma\right)$, so the separation of the eigenvalues with respect to the unit circle is particularly unfavorable, leading to a small-scale, severely ill-conditioned problem. 
\item[Example 4b] The same as Example 4a, but with $\sigma = 10^{-5}$. This choice produces a small-scale, mildly ill-conditioned problem.
\end{description}

For each test, we report the relative residual in the Euclidean norm
\[
\mathrm{Res} = \frac{\norm{DX +X^\top A - X^\top BX + C}}{\norm{D}\norm{X} + \norm{X}\norm{A} + \norm{X}\norm{B}\norm{X} + \norm{C}},
\]
and for the smaller examples also the forward relative error
\[
\mathrm{Err} = \frac{\norm{X - X_{\mathrm{exact}}}}{\norm{X_{\mathrm{exact}}}}
\]
with respect to a reference solution computed with the Doubling Algorithm, run in higher precision using Matlab's \verb!vpa! command. This allows to test both the backward and the forward stability of the algorithms.

The results are given in Tables~\ref{fig:ex1}--\ref{fig:ex4b}.
\begin{figure}
    \centering
    \begin{tabular}{llll}
\toprule
& \textbf{Res} & \textbf{Err} & \textbf{Notes} \\
\midrule
\textbf{QZ} & 7.973374e-16 & 4.734635e-15 &  \\ 
\textbf{PalQZ} & 7.051521e-16 & 2.190775e-15 &  \\ 
\textbf{DA} & 8.098123e-17 & 1.735516e-16 & $k = 8$ \\ 
\textbf{CR} & 5.547051e-17 & 1.969563e-16 & $k = 8$ \\ 
\textbf{CR2} & 1.042184e-15 & 9.124629e-15 & $k = 8$ \\ 
\textbf{PDA} & 5.785320e-16 & 1.026009e-15 & $k = 9$ \\ 
\textbf{Int} & 4.568426e-16 & 1.107647e-15 & $k = 8$ \\ 
\bottomrule
    \end{tabular}
    \caption{Numerical results on Example 1.}
    \label{fig:ex1}
\end{figure}
\begin{figure}
    \centering
    \begin{tabular}{lll}
\toprule
& \textbf{Res} & \textbf{Notes} \\
\midrule
\textbf{QZ} & 2.426140e-15 &  \\ 
\textbf{PalQZ} & 4.181695e-15 &  \\ 
\textbf{DA} & 1.210196e-16 & $k = 9$ \\ 
\textbf{CR} & 2.820211e-16 & $k = 9$ \\ 
\textbf{CR2} & 1.483258e-15 & $k = 9$ \\ 
\textbf{PDA} & 8.703012e-16 & $k = 9$ \\ 
\textbf{Int} & 5.529665e-16 & $k = 8$ \\ 
\bottomrule
    \end{tabular}
    \caption{Numerical results on Example 2a.}
    \label{fig:ex2a}
\end{figure}
\begin{figure}
    \centering
    \begin{tabular}{lll}
\toprule
& \textbf{Res} & \textbf{Notes} \\
\midrule
\textbf{QZ} & 7.412436e-15 &  \\ 
\textbf{PalQZ} & 5.335863e-15 &  \\ 
\textbf{DA} & 4.133160e-15 & $k = 11$ \\ 
\textbf{CR} & 4.417749e-15 & $k = 11$ \\ 
\textbf{CR2} & 3.048525e-13 & $k = 11$ \\ 
\textbf{PDA} & 5.856347e-15 & $k = 11$ \\ 
\textbf{Int} & 6.607131e-16 & $k = 11$\\ 
\bottomrule
    \end{tabular}
    \caption{Numerical results on Example 2b.}
    \label{fig:ex2b}
\end{figure}
\begin{figure}
    \centering
    \begin{tabular}{llll}
\toprule
& \textbf{Res} & \textbf{Err} & \textbf{Notes} \\
\midrule
\textbf{QZ} & 4.662032e-08 & 2.238634e-06 &  \\ 
\textbf{PalQZ} & 2.109338e-17 & 6.571568e-15 &  \\ 
\textbf{DA} & 1.162783e-07 & 5.583500e-06 & $k = 38$ \\ 
\textbf{CR} & 5.001176e-08 & 2.401486e-06 & $k = 38$ \\ 
\textbf{CR2} & --- & --- & Did not converge \\ 
\textbf{PDA} & 1.206223e-07 & 5.812851e-06 & $k = 38$ \\ 
\textbf{Int} & --- & --- & Would take too long to run with $k=38$ \\ 
\bottomrule
    \end{tabular}
    \caption{Numerical results on Example 4a.}
    \label{fig:ex4a}
\end{figure}
\begin{figure}
    \centering
    \begin{tabular}{llll}
\toprule
& \textbf{Res} & \textbf{Err} & \textbf{Notes} \\
\midrule
\textbf{QZ} & 5.004467e-12 & 1.488394e-09 &  \\ 
\textbf{PalQZ} & 4.759728e-17 & 6.526349e-15 &  \\ 
\textbf{DA} & 4.328464e-12 & 1.288456e-09 & $k = 22$ \\ 
\textbf{CR} & 7.349975e-11 & 2.187141e-08 & $k = 22$ \\ 
\textbf{CR2} & 1.337969e-03 & 5.442744e-01 & $k = 23$ \\ 
\textbf{PDA} & 1.612961e-12 & 4.823130e-10 & $k = 24$ \\ 
\textbf{Int} & 2.132475e-14 & 2.808352e-11 & $k = 22$ \\ 
\bottomrule
    \end{tabular}
    \caption{Numerical results on Example 4b.}
    \label{fig:ex4b}
\end{figure}
One sees that the accuracy obtained with the palindromic QZ algorithm is unmatched by the other algorithms. Most of the iterative algorithms obtain similar accuracy results to the (unstructured) QZ algorithm, despite some of them formally exploiting symmetry; in particular, it is interesting to note that the algorithm based on the integral computation with the trapezoidal rule is slightly more accurate, although significantly slower: its cost scales as $\mathcal{O}(2^k n^3)$ rather than $\mathcal{O}(kn^3)$.

The negative outlier in terms of accuracy is algorithm \textbf{CR2}, which is the one based on the symmetric quadraticization~\eqref{eq:infs}. Indeed, one sees that this infinite system of equations is obtained with a sort of infinite-dimensional analogue of the normal equations 
\begin{multline*}
\begin{bmatrix}
\mathcal{M}^\top \mathcal{M} & -\mathcal{M}^{2\top}\\
-\mathcal{M}^{2} & \mathcal{M}^\top \mathcal{M} + \mathcal{M} \mathcal{M}^\top & -\mathcal{M}^{2\top}\\
 &   -\mathcal{M}^{2} & \mathcal{M}^\top \mathcal{M} + \mathcal{M} \mathcal{M}^\top & \ddots\\
 && \ddots & \ddots 
\end{bmatrix}
\\=
\begin{bmatrix}
\mathcal{M}^\top\\
-\mathcal{M} & \mathcal{M}^\top\\
& \ddots & \ddots
\end{bmatrix}
\begin{bmatrix}
\mathcal{M} & -\mathcal{M}^\top \\
& \mathcal{M} & -\mathcal{M}^\top \\
& & \ddots & \ddots \\
\end{bmatrix}.
\end{multline*}
Hence a lower accuracy on ill-conditioned examples is expected, analogously to what happens in other algorithms based on normal equations, which notoriously may exhibit a ``squaring the condition number'' phenomenon; see, e.g.,~\cite[Section~5.3.2]{gvl} for a discussion.

Finally, a remark on the number of iterations: $k\approx$ 8--10 is sufficient for most examples without problematic eigenvalues, while the number of iterations may increase greatly if there are eigenvalues very close to the unit circle. This slowdown is expected from the theory, as the convergence speed of doubling algorithms is $\mathcal{O}(\sigma^{2^k})$.

\subsection{Timing of the swap procedure}

We test the swap procedure in Section \ref{sec:sw} and compare it with the swap procedure in \cite{bimp}. The theoretical analysis shows that the proposed technique requires $\mathcal{O}(n^3)$ ops in the worst case, while the technique of \cite{bimp} requires $\mathcal{O}(n^4)$ ops in the worst case. The worst case is the same for both methods and corresponds to an antitriangular pencil where all the eigenvalues smaller than one are in the top right part. 

We test the algorithms on random antitriangular matrix pencils $R+zR^\top $, where $R$ is obtained with the MATLAB command \verb+flipud(triu(randn(n)))+, where $n$ is the size. 
In Table \ref{tab:swap2}, we report the times measured for our algorithm to show how the predicted cost is realized in practice on an average case. The reported times for each $n$ are obtained by averaging $5$ runs of the algorithm, and the average number of swaps is provided.
The results of the algorithm in \cite{bimp} are not reported since it fails to reorder the antitriangular form for matrices of size greater than about $100\times 100$.

\begin{table}
\begin{center}
\begin{tabular}{ccc}
\hline
\textbf{n} & \textbf{Time (s)} & \textbf{Swaps}\\\hline
32 & 0.0021 & 54.2 \\
64 & 0.0082 & 282.2 \\
128 & 0.023 & 1094.2 \\
256 & 0.089 & 3901.4\\
512 & 0.489 & 16230.0\\
1024 & 3.442 & 66202.2\\
2048 & 30.13 & 263385.2\\
4096 & 290.1 & 1032142.0\\
\hline
\end{tabular}
\caption{Results of the experiment on the reordering  algorithm for an antitriangular pencil.}
\label{tab:swap2}
\end{center}
\end{table}

\section{Conclusions}\label{sec:con}

In this paper we obtain both theoretical and computational results on $\top$-NARE, relying on their deflating subspace formulation. In the experiments, we compare a selection of various new doubling-type algorithms to the existing algorithms based on deflating subspace computations. Even if these variants take advantage of the symmetry of the problem computationally, the numerical results show that the accuracy of the palindromic QZ algorithm is unparalleled. These experiments also highlight the benefits of a faster swapping strategy to complete the algorithm, which is one of our results. Finding a doubling variant that fully exploits the structure and matches the accuracy of the palindromic QZ remains an open research problem.

\section*{Acknowledgments} This work was partly supported by INdAM (Istituto Nazionale di Alta Matematica) through a GNCS project, and by the University of Pisa's project PRA\_2020\_61. Some of the software used for the numerical experiments involving the \textbf{PalQZ} method is originally due to D.~Palitta, F.~Dopico and C.~Schroeder; we are grateful to these authors.

\vspace{1cm}

\bibliographystyle{abbrv}
	\bibliography{tnare}

\begin{thebibliography}{10}

\bibitem{aishima}
K.~Aishima.
\newblock A quadratically convergent algorithm based on matrix equations for
  inverse eigenvalue problems.
\newblock {\em Linear Algebra and its Applications}, 542:310--333, 2018.
\newblock Proceedings of the 20th ILAS Conference, Leuven, Belgium 2016.

\bibitem{BaiDG97}
Z.~Bai, J.~Demmel, and M.~Gu.
\newblock An inverse free parallel spectral divide and conquer algorithm for
  nonsymmetric eigenproblems.
\newblock {\em Numer. Math.}, 76(3):279--308, 1997.

\bibitem{bimp}
P.~Benner, B.~Iannazzo, B.~Meini, and D.~Palitta.
\newblock Palindromic linearization and numerical solution of nonsymmetric
  algebraic {T}-{Riccati} equations.
\newblock {\em BIT Numerical Mathematics}, 62(4):1649 – 1672, 2022.

\bibitem{bp20}
P.~Benner and D.~Palitta.
\newblock {On the solution of the nonsymmetric T-Riccati equation}.
\newblock {\em Electron. Trans. Numer. Anal.}, 54:68--88, 2021.

\bibitem{bv}
M.~Benzi and M.~Viviani.
\newblock Solving cubic matrix equations arising in conservative dynamics.
\newblock {\em Vietnam J. Math.}, 2022.

\bibitem{blm}
D.~A. Bini, G.~Latouche, and B.~Meini.
\newblock {\em Numerical methods for structured {M}arkov chains}.
\newblock Numerical Mathematics and Scientific Computation. Oxford University
  Press, New York, 2005.
\newblock Oxford Science Publications.

\bibitem{bm:cr}
D.~A. Bini and B.~Meini.
\newblock The cyclic reduction algorithm: from {P}oisson equation to stochastic
  processes and beyond. {I}n memoriam of {G}ene {H}. {G}olub.
\newblock {\em Numer. Algorithms}, 51(1):23--60, 2009.

\bibitem{borobia}
A.~Borobia, R.~Canogar, and F.~De~Ter\'an.
\newblock On the consistency of the matrix equation {$X^\top AX= B$} when {$B$}
  is symmetric.
\newblock {\em Mediterranean J. Math.}, 18(2), 2021.

\bibitem{borobia2}
A.~Borobia, R.~Canogar, and F.~De~Ter\'an.
\newblock The equation {$X^T AX=B$} with {$B$} skew-symmetric: how much of a
  bilinear form is skew-symmetric?
\newblock {\em Linear and Multilinear Algebra}, 2022.

\bibitem{bk}
R.~Byers and D.~Kressner.
\newblock Structured condition numbers for invariant subspaces.
\newblock {\em SIAM J. Matrix Anal. Appl.}, 28(2):326--347, 2006.

\bibitem{dipr19}
F.~De~Ter\'{a}n, B.~Iannazzo, F.~Poloni, and L.~Robol.
\newblock Nonsingular systems of generalized {S}ylvester equations: an
  algorithmic approach.
\newblock {\em Numer. Linear Algebra Appl.}, 26(5):e2261, 29, 2019.

\bibitem{edb}
C.~Estatico and F.~Di~Benedetto.
\newblock Shift-invariant approximations of structured shift-variant blurring
  matrices.
\newblock {\em Numerical Algorithms}, 62(4):615--635, 2013.

\bibitem{gant}
F.~R. Gantmacher.
\newblock {\em The theory of matrices. {V}ols. 1, 2}.
\newblock Chelsea Publishing Co., New York, 1959.
\newblock Translated by K. A. Hirsch.

\bibitem{gardiner}
J.~D. Gardiner and A.~J. Laub.
\newblock A generalization of the matrix-sign-function solution for algebraic
  {R}iccati equations.
\newblock {\em International Journal of Control}, 44(3):823--832, 1986.

\bibitem{glr}
I.~Gohberg, P.~Lancaster, and L.~Rodman.
\newblock {\em Matrix Polynomials}.
\newblock SIAM, 2009.

\bibitem{gvl}
G.~H. Golub and C.~F. van Loan.
\newblock {\em Matrix Computations}.
\newblock John Hopkins University Press, Baltimore, fourth edition, 2013.

\bibitem{highamfun}
N.~J. Higham.
\newblock {\em Functions of matrices}.
\newblock Society for Industrial and Applied Mathematics (SIAM), Philadelphia,
  PA, 2008.
\newblock Theory and computation.

\bibitem{ikr}
K.~Ikramov.
\newblock On the solvability of a certain class of quadratic matrix equations.
\newblock {\em Doklady Mathematics}, 89(2):162 – 164, 2014.

\bibitem{ksw09}
D.~Kressner, C.~Schr{\"o}der, and D.~S. Watkins.
\newblock Implicit {QR} algorithms for palindromic and even eigenvalue
  problems.
\newblock {\em Numer. Algorithms}, 51(2):209--238, 2009.

\bibitem{lancasterrodman-are}
P.~Lancaster and L.~Rodman.
\newblock {\em Algebraic {R}iccati equations}.
\newblock Oxford Science Publications. The Clarendon Press, Oxford University
  Press, New York, 1995.

\bibitem{pda}
T.~Li, C.-Y. Chiang, E.~K.-w. Chu, and W.-W. Lin.
\newblock The palindromic generalized eigenvalue problem {$A^*x=\lambda Ax$}:
  numerical solution and applications.
\newblock {\em Linear Algebra Appl.}, 434(11):2269--2284, 2011.

\bibitem{m4}
D.~S. Mackey, N.~Mackey, C.~Mehl, and V.~Mehrmann.
\newblock Numerical methods for palindromic eigenvalue problems: computing the
  anti-triangular {S}chur form.
\newblock {\em Numer. Linear Algebra Appl.}, 16(1):63--86, 2009.

\bibitem{bea:extr}
B.~Meini.
\newblock Efficient computation of the extreme solutions of {$X+A^*X^{-1}A=Q$}
  and {$X-A^*X^{-1}A=Q$}.
\newblock {\em Math. Comp.}, 71(239):1189--1204, 2002.

\bibitem{miy}
S.~Miyajima.
\newblock Fast enclosure for the minimal nonnegative solution to the
  nonsymmetric {$T$}-{Riccati} equation.
\newblock {\em Calcolo}, 59(3), 2022.

\bibitem{feast}
E.~Polizzi.
\newblock Density-matrix-based algorithm for solving eigenvalue problems.
\newblock {\em Physical Review B}, 79(11):115112, 2009.

\bibitem{stewart72}
G.~W. Stewart.
\newblock On the sensitivity of the eigenvalue problem $ax = \lambda bx$.
\newblock {\em SIAM Journal on Numerical Analysis}, 9(4):669--686, 1972.

\bibitem{vor}
Y.~O. Vorontsov.
\newblock Solvability conditions for the matrix equation {$X^TDX +AX + X^TB + C
  = 0$}.
\newblock {\em Computational Mathematics and Mathematical Physics}, 55(4):546
  – 548, 2015.

\bibitem{vor2}
Y.~O. Vorontsov and K.~D. Ikramov.
\newblock Numerical algorithm for solving quadratic matrix equations of a
  certain class.
\newblock {\em Computational Mathematics and Mathematical Physics}, 54(11):1643
  – 1646, 2014.

\bibitem{yuan}
Y.~Yuan, Y.~Chen, and L.~Zhou.
\newblock Solutions of a second-order conjugate matrix equation.
\newblock {\em Linear and Multilinear Algebra}, 69(9):1645 – 1656, 2021.

\bibitem{yuan1}
Y.~Yuan, L.~Liu, H.~Zhang, and H.~Liu.
\newblock The solutions to the quadratic matrix equation {$X^*AX+B^*X+D=0$}.
\newblock {\em Applied Mathematics and Computation}, 410, 2021.

\end{thebibliography}

\end{document}